%%%%%%%%%%%%%%%%%%%%%%%%%%%%%%%%%%%%%%%%%%%%%%%%    
%
%        THIS IS A  PLAIN TeX FILE
%
%%%%%%%%%%%%%%%%%%%%%%%%%%%%%%%%%%%%%%%%%%%%%%%%

\magnification=1200

\font\titfont=cmr10 at 12 pt

\font\headfont=cmr10 at 12 pt

%\font\AAA=Times.dfont  at 12pt
 %\font\BBB=Times.dfont at 8pt

%\font\AAA=cmr10 at 12pt
%\font\BBB=cmr10 at 8pt

\overfullrule=0in

\def\boxit#1{\hbox{\vrule
 \vtop{%
  \vbox{\hrule\kern 2pt %
     \hbox{\kern 2pt #1\kern 2pt}}%
   \kern 2pt \hrule }%
  \vrule}}

  \def\harr#1#2{\ \smash{\mathop{\hbox to .3in{\rightarrowfill}}\limits^{\scriptstyle#1}_{\scriptstyle#2}}\ }

 \def\GG{{{\bf G} \!\!\!\! {\rm l}}\ }

\def\bll{I \!\! L}

\def\bra#1#2{\langle #1, #2\rangle}
\def\bbf{{\bf F}}
\def\bbj{{\bf J}}
\def\Jtn{{\bbj}^2_n}  \def\JtN{{\bbj}^2_N}  \def\JoN{{\bbj}^1_N}
\def\jt{j^2}
\def\jtx{\jt_x}
\def\Jt{J^2}

\def\ss{\subset}

\def\half{\hbox{${1\over 2}$}}
\def\smfrac#1#2{\hbox{${#1\over #2}$}}

\def\dim{{\rm dim}}
\def\dist{{\rm dist}}
\def\codim{{\rm codim}}

\def\rank{{\rm rank}}

\def\Hess{{\rm Hess}}

\def\trace{{\rm trace}}
\def\tr{{\rm tr}}
\def\max{{\rm max}}

\def\span{{\rm span\,}}
\def\Hom{{\rm Hom\,}}
\def\det{{\rm det}}
\def\End{{\rm End}}
\def\Sym{{\rm Sym}^2}

\def\arr{\longrightarrow}

\def\rn{\bbr^n}

\def\Int{{\rm Int}}

\def\Symn{{\Sym(\rn)}}

\def\Theorem#1{\medskip\noindent {\bf THEOREM \bf #1.}}
\def\Prop#1{\medskip\noindent {\bf Proposition #1.}}
\def\Cor#1{\medskip\noindent {\bf Corollary #1.}}
\def\Lemma#1{\medskip\noindent {\bf Lemma #1.}}
\def\Remark#1{\medskip\noindent {\bf Remark #1.}}
\def\Note#1{\medskip\noindent {\bf Note #1.}}
\def\Def#1{\medskip\noindent {\bf Definition #1.}}

\def\Ex#1{\medskip\noindent {\bf Example \bf    #1.}}

\def\pf{\medskip\noindent {\bf Proof.}\ }
\def\qed{\hfill  $\vrule width5pt height5pt depth0pt$}

   \def\cp{{\cal P}}

   \def\cn{{\cal N}}

\def\cp{{\cal P}}
\def\cf{{\cal F}}

\def\vf{\varphi}

\def\wt{\widetilde}

\def\and{\qquad {\rm and} \qquad}
\def\arr{\longrightarrow}

\def\bbr{{\bf R}}\def\bbh{{\bf H}}
\def\bbc{{\bf C}}

\def\a{\alpha}
\def\b{\beta}
\def\d{\delta}
\def\e{\epsilon}

\def\g{\gamma}

\def\l{\lambda}

\def\D{\Delta}
\def\L{\Lambda}
\def\G{\Gamma}
\def\O{\Omega}

\def\Symn{\Sym(\rn)}
 
\def\USC{{\rm USC}}
\def\fa{{\rm\ \  for\ all\ }}

\def\AA{1}
\def\BB{2}
\def\CC{3}
\def\DD{4}
\def\EE{4}
\def\FF{5}
\def\HH{6}

\def\JJ{8}
\def\LL{7}
 \def\KK{A}
\def\MM{8}
\def\NN{9}

\centerline{\titfont  THE  RESTRICTION THEOREM }
\smallskip

\centerline{\titfont  FOR FULLY NONLINEAR SUBEQUATIONS}
\bigskip

\centerline{\titfont F. Reese Harvey and H. Blaine Lawson, Jr.$^*$}
\vglue .9cm
\smallbreak\footnote{}{ $ {} \sp{ *}{\rm Partially}$  supported by
the N.S.F. } 

\vskip .2in
 
\centerline{\bf ABSTRACT} \medskip
  \font\abstractfont=cmr10 at 10 pt
{{\parindent= .43in
\narrower\abstractfont \noindent

We address the restriction problem for viscosity subsolutions of a fully nonlinear
PDE on a manifold $Z$. The constraints on the
restrictions of smooth subsolutions to a submanifold $X\ss Z$
 determine a restricted subequation on $X$.
The problem is to show that general (upper semi-continuous) subsolutions restrict to satisfy the same
constraints.  We first prove an elementary result which, in theory, can be applied to
any subequation.  Then two definitive results are obtained. The first applies to any
``geometrically defined'' subequation, and the second to  any subequation which 
can be transformed to a constant coefficient (i.e., euclidean) model.  This provides a long
list of geometrically and analytically interesting cases where restriction holds.

}}

\vskip.5in

\centerline{\bf TABLE OF CONTENTS} \bigskip

{{\parindent= .1in\narrower\abstractfont \noindent

\qquad 1. Introduction.\smallskip

\qquad \BB.     Nonlinear Potential Theory    \smallskip

\qquad \CC.    Introduction to Restriction -- The Geometric Case in $\rn$.  \smallskip

\qquad \EE.     The General Restriction Theorem.  \smallskip

\qquad \FF.   First Applications.\smallskip

\qquad \HH.   The Geometric  Restriction Theorem.\smallskip

%\qquad \II.   Geometric Examples on Riemannian Manifolds.\smallskip

%\qquad \JJ.   $\GG$-flat Submanifolds.\smallskip

\qquad \LL.   Jet Equivalence of Subequations \smallskip

\qquad \MM.   The  Restriction Theorem for Subequations Derivable from a Euclidean 

\qquad \ \ \ \ \ Model.\smallskip

\qquad \NN.   Applications of this Last  Restriction Theorem.\smallskip

}}

\vskip .1in

{{\parindent= .3in\narrower

\noindent
{\bf Appendices: }\medskip

A.  Elementary  Examples Where Restriction Fails
\smallskip

B. Restriction of Sets of Quadratic Forms.

\smallskip
 
C. Extension Results.

\smallskip

}}

\vfill\eject

\noindent{\headfont \AA.\  Introduction}
\medskip

This paper is concerned with the  restrictions of   subsolutions
of a fully nonlinear elliptic partial differential equation to  submanifolds.
In most cases this topic is uninteresting because the   restricted
functions satisfy no constraints. Moreover, even when there are 
 constraints, this will occur only on certain submanifolds.
 Nonetheless,  there are cases, in fact many cases,  where the restriction question
is quite interesting.  Important classical examples  are
the plurisubharmonic functions in several complex variable theory,
and their analogues in calibrated geometry.
 The principle aim of this paper is to  study
the foundations of the restriction problem.
We prove a general Restriction Theorem which applies to all cases, but  
whose ``restriction hypothesis'' must be verified.  We then obtain definitive results
in two general situations, each followed with a  series of applications. 
First, if the constraints are ``determined geometrically'', the applications
 come from potential theory developed in calibrated and other geometries 
(cf. [HL$_{2,3}$]).  In the second situation the constraints are locally derivable 
from a constant coefficient (euclidean) model.  Here the applications
come from universal subequations in riemannian geometry
  (cf. [HL$_{6,7}$]).  Yet another application  will be to the study of the intrinsic 
   potential theory on almost complex manifolds (without use of a hermitian metric) [HL$_8$].

We begin with a note about our approach to this problem.
 Traditionally, a second-order partial differential equation (or subequation)  is a constraint
on the full second derivative (or 2-jet) of a function $u$
 imposed by using a function  $f(x,u, Du, D^2u)$ and setting  $f=0$ (or $f \geq0$).  
We have found it more enlightening to work directly with the  subsets   of the 
2-jet space corresponding to these conditions (cf.[K]), and we have systematically explored
this viewpoint in recent papers [HL$_{4,5,6}$].  
(A succinct comparison of  our subset  approach  with the standard one is given in 
a Pocket Dictionary in  [HL$_9$, App. A].)
This geometric
formulation is often more natural and   has several  distinct
advantages. To begin, it makes the equation completely canonical.
It  clarifies a number of  classical conditions, such as the condition of
degenerate ellipticity. It underlines an inherent duality in the subject, which
in turn clarifies the necessary boundary geometry for solving the Dirichlet
problem.

It also simplifies and clarifies certain natural operations, in particular those of restriction
and addition.  

To be more concrete,  let's begin with  a closed subset $F$ of the space of 2-jets over a domain
$Z \ss\rn$, which we assume to satisfy the very weak ellipticity condition (\BB.4) below, called {\sl positivity}.
Such a set will be called a {\sl subequation}.
Then a function $u \in C^2(Z )$ is called {\sl $F$-subharmonic} if its 2-jet $J^2_xu \in F$
 for all $x$. This concept can be extended to upper semi-continuous functions
$u$ by using the following viscosity approach (cf.[CIL], [C]).  
We say that a function $\vf$ which is $C^2$ near $x\in Z$  is   {\sl a  test function
for $u$ at $x$} if $u\leq \vf$ near $x$ and $u(x)=\vf(x)$.
Then a function $u\in \USC(Z )$ is {\sl $F$-subharmonic} if for 
 each  test function $\vf$ for $u$ at  any $x\in Z$, one has $J^2_x\vf\in F$.

Suppose now that
$$
F\ss J^2(Z)
$$ 
 is a subequation and $i:X\ss Z$ is a submanifold of  $Z$. 
  Then there is a naturally induced subequation
$$
H \equiv \overline{i^*F} \ss J^2(X)
$$
where $i^*F$ is given by the restriction of 2-jets (which is induced by the  restriction of smooth functions).
  By definition it has the property that 
 for $\vf\in C^2(Z)$
 $$
 \vf \ \ {\rm is \ } F{\rm -subharmonic}
 \qquad\Rightarrow\qquad
 \vf\bigr|_X \ \ {\rm is \ }    \overline{i^*F}{\rm -subharmonic}
 \eqno{(\AA.1)}
 $$
 
As mentioned before, for general  $F$ and $X$ the induced subequation
is uninteresting.  This is because generically $i^*F=J^2(X)$, and so no constraints
are placed on restrictions of $F$-subharmonic functions. This leads to two natural problems.

\medskip
\noindent
{\bf Problem 1:}   \  Identify non-vacuous cases and calculate the induced subequation $\overline{i^*F}$.
\medskip

Frequently $i^*F$ is closed,  but not always (see Examples \FF.5 and B.6).

Once Problem 1 is accomplished, we have the second, more difficult task of determining
whether restriction holds.
\medskip
\noindent
{\bf Problem 2:}   \  Find conditions under which the restriction statement (\AA.1)  
extends to upper semi-continuous functions.
\medskip

In the classical case coming from several complex variable theory, the subequation
$F$ is defined by requiring that the complex hermitian part of the hessian matrix be 
non-negative. 
Here the most interesting submanifolds are the complex curves, and in this case
the restricted subequation is the conformal Laplacian.  Thus the prototype of our main result is
the theorem which says that a function which is  plurisubharmonic in the viscosity sense
is the same as a function whose restriction to every complex curve is subharmonic.
In fact,  the corresponding statement has recently been established 
for almost complex manifolds by using one of our main results Theorem \MM.1.  This application 
is presented in a separate paper [HL$_8$].
(See Note 1.1 below.)

An even more basic case is the real analogue, which states that a function is convex
in the viscosity sense if and only if its restriction to each affine line is convex.
 
 These classical cases extend to branches of the  homogeneous Monge-Amp\`ere
 equation and the concept of $q$-convexity. The whole story carries over
 to the important complex and quaternionic settings where much work has been done.
See Note 1.2  for a more detailed discussion of this and some generalizations.

There are many  other general cases in which the outcome of Problem 1 is known and interesting.
Some come from potential theory developed in calibrated and other geometries 
(cf. [HL$_{2,3}$]).  Others come from universal subequations in riemannian geometry
 and on manifolds with topological $G$-structures (cf. [HL$_{6,7}$]). These will all be investigated here.

 We begin the paper   with definitions and a brief review of potential theory for 
 fully nonlinear subequations.  In order to introduce and motivate the restriction
 problem,  we first examine it for ``geometrically determined  subequations''
  in euclidean space.  These are subequations
 $F_\GG$ determined by the condition $\trace\{D^2 u\bigr|_W\}\geq 0$ for all 
 $p$-planes $W$ in a given fixed subset $\GG\ss G(p,\rn)$ of the grassmannian 
 of $p$-planes in $\rn$. This, of course,  includes the classical case of plurisubharmonic functions
 in complex analysis where $\GG \equiv G_\bbc(1,\bbc^n) \ss G(2,\bbr^{2n})$.
 
 In Section \EE \ we  prove   a basic  elementary theorem. 
 For a given subequation $F\ss J^2(Z)$ and submanifold $i:X\ss Z$, we formulate
 a {\sl Restriction Hypothesis} and prove the following.

\smallskip
\noindent
{\bf The Restriction Theorem \EE.2.}  {\sl  Suppose $u\in \USC(Z)$.  
Assume that $F$ satisfies  the restriction hypothesis.  Then}
 $$
 u\in F(Z)  \qquad \Rightarrow\qquad u\bigr|_X \in  (\overline{i^*F})(X).
 $$

The proof parallels a proof of Crandall [C, Lemma 4.1].

 This result is then applied throughout the rest of the paper.
 
 In Section \FF\  we make some immediate but important applications.    
 Two of them are prototypes for  the main restriction theorems in this paper.
The first presented is the following.
{\sl  Suppose $F\ss J^2(Z)$ is a translation-invariant, i.e., constant coefficient, 
 subequation on an open set $Z\ss \bbr^n$. Then
for $u\in \USC(Z)$, 
\medskip
\centerline
{
$u$ is $F$-subharmonic on $Z$ \qquad$\Rightarrow\qquad u\bigr|_X$ is \ 
$\overline { i^*F}$-subharmonic on $X$.
}
}
\medskip
\noindent
For the second prototype we  establish restriction for $\GG$-plurisubharmonic functions, 
to affine $\GG$-planes (defined below).
In addition to these two prototypes we establish restriction for 
general linear subequations under a (necessary)
{\sl linear restriction hypothesis}.  This result becomes important in later applications.
Finally, we examine restriction for first-order equations.

 In Section \HH  \ we establish our quite general and definitive Restriction Theorem
 \HH.6.   A special case is the following.
 Let $Z$ be a riemannian manifold of dimension $n$ and  $\GG\ss G(p, TZ)$ 
a closed subset of the bundle of tangent $p$-planes on $Z$. 
Assume   that $\GG\ss G(p, TZ)$ 
admits a smooth neighborhood retraction which preserves the fibres 
of the projection $G(p, TZ) \to Z$.
Then $\GG$ determines a 
natural subequation $F$ on $Z$ defined by the condition that
$$
\trace\left\{\Hess \, u\bigr|_W \right\} \ \geq \ 0 \fa W \in \GG.
$$
where $\Hess\, u$ denotes the riemannian hessian of $u$.
(See (\NN.3) and  [HL$_2$] for examples and details.) 
The corresponding  $F$-subharmonic functions
are again called $\GG${\bf -plurisubhamronic functions}.

A {\bf $\GG$-submanifold} of $Z$ is defined to be a $p$-dimensional submanifold
$X\ss Z$ such that $T_xX \in \GG$ for all $x\in X$.

\Theorem{\HH.4} {\sl Let $X\ss Z$ be a $\GG$-submanifold which is minimal (mean curvature zero).
Then restriction to $X$ holds for  $F$.  In other words, the restriction of any
 $\GG$-plurisubharmonic function to $X$ is subharmonic in the induced riemannian  metric
 on $X$.}
 \medskip

In the general result, Theorem \HH.6, the  submanifold is allowed to have  dimension $>p$.

 In Sections \LL \ and \MM \ we formulate a quite different restriction result, based on the 
idea of jet equivalence. The notion of jet  equivalence of subequations was introduced in [HL$_6$, \S 4]
where it greatly extended the applicability of basic results. This notion is recalled 
in   Section \LL \ and then refined to the relative case.
 We then prove the following for an open subset $Z\ss \bbr^N$ containing an embedded submanifold
 $i:X\hookrightarrow Z$.

\Theorem{\MM.1}  {\sl Suppose that $F\ss J^2(Z)$ 
a subequation.  Assume that 
  $F$ is  locally jet equivalent  modulo $X$ to a constant coefficient
subequation $\bbf$.  Then 
$H\equiv \overline{i_X^*F}$ is locally jet equivalent to the constant coefficient subequation 
$\bbh \equiv \overline{i^* \bbf}$. Moreover, restriction holds.  That is, }
$$
u \ \ {\sl is}\ F\ {\sl subharmonic\ on\ }\ Z\qquad\Rightarrow\qquad 
u\bigr|_X \ \ {\sl is}\ H\ {\sl subharmonic\ on\ }\  X
$$

This theorem has a number of interesting applications.  One is the following.

\Theorem {\MM.2} {\sl
Let $Z$ be a riemannian manifold of dimension $N$ and $F\ss J^2(Z)$ a subequation canonically
determined by an O$_N$-invariant universal subequation $\bbf\ss\bbj^2_N$ (see \S \MM).
Then restriction  holds for $F$ on any totally geodesic submanifold $X\ss Z$.
}
\medskip

This result extends to subequations defined by $G$-invariant 
subsets of $\bbj^2_N = \bbr\times \bbr^N\times \Sym(\bbr^N)$ on manifolds
with topological $G$-structure.

In constrast to Theorem \HH.6, a riemannian metric is {\sl not} required 
in Theorem \MM.1, so that it can be applied as follows.

\medskip
\noindent
{\bf Note 1.1. (Almost Complex Manifolds and the Pali Conjecture).} 
Another application of  Theorem \MM.1 is to the study of potential
theory on almost complex manifolds in the absence of any hermitian metric.
In this case there is still an intrinsically defined subequation,
but it is not geometrically defined in the sense of Section \HH.
The corresponding subharmonic functions are proved in [HL$_8$] (Theorem 6.2)
to be exactly those upper semi-continuous functions 
whose restrictions to complex curves are subharmonic.
This  is then used to establish the full version  of a conjecture of Nefton Pali [P].
(See Theorem 8.2 of [HL$_8$].)
The restriction theorem is central to this work.

\medskip
\noindent
{\bf Note 1.2. (Branches of the Homogeneous Monge-Amp\`ere Equation and $q$-Convexity).} 
 The  classical cases of convex and plurisubharmonic functions discussed above
 can be extended as follows.
For $A\in\Symn$, let $\l_1(A) \leq \l_2(A) \leq
 \cdots \leq \l_n(A)$  denote its ordered eigenvalues. Then the $q^{\rm th}$ {\sl branch} of 
 $\det(D^2 u)=0$ is the equation $\l_q(D^2 u)=0$. Its associated subequation
 $\L_q \equiv \{A : \l_q(A) \geq0\}$ is the condition of $q$-{\sl convexity}.
The u.s.c. $\L_q$-subharmonic functions will be called  $q$-{\sl convex}.
When $q=1$ these are just the convex functions, and the first branch of the
homogeneous Monge-Amp\`ere Equation is the classical one treated by 
Alexandrov [Al].  When $q=n$ the  n-convex functions
 are the  {\sl subaffine} functions introduced in [HL$_4$]. (See Def. \BB.6 and Prop. \BB.7).
 For general $q$ our restriction results apply to prove that for an open set 
 $X\ss\rn$:
 
 \noindent
 {\bf Theorem \FF.3.}
  {\sl A function $u\in \USC(X)$ is $q$-convex if and only if its restriction to 
 every affine $q$-plane is subaffine.}

 This entire story carries over to the complex case in $\bbc^n = \bbr^{2n}$ 
  by replacing $A$ with its hermitian symmetric part $A_\bbc \equiv \half(A- JAJ)$.
  %where $J$ denotes multiplication by $i$ on $\bbr^{2n}$. 
Here one studies branches of the homogeneous {\sl complex} Monge-Amp\`ere equation. 
The first branch is the classical one (cf. {B},[BT]) and the 1-convex functions are the
plurisubharmonic functions discussed above.  At the other end, the largest branch
of $n$-convex functions can be characterized as 
those which are ``sub-the-pluriharmonics'' (see Def. \FF.13 and Prop. \FF.14).
The case of general $q$ has received 
much attention in complex analysis (e.g.   [HM], [S]). 
Our restriction results show that for an open set 
 $X\ss\bbc^n$:
 
 \noindent
 {\bf Theorem \FF.16}
  {\sl A function $u\in \USC(X)$ is $q$-convex (in the complex sense)  if and only if its restriction to 
 every complex affine $q$-plane is sub-the-pluriharmonics.}
  
 There is also   an interesting quaternionic analogue  
 of the Monge-Amp\`ere equations (cf. [A$_*$], [AV])
and associated   $q$-convex functions  
 [HL$_{2,4,6}$]. The assertions above generalize to this case.
Results for inhomogeneous equations on manifolds are
given in  Example \NN.7.

\medskip
\noindent
{\bf Remark 1.3.} The Restriction Hypothesis discussed in Section 4 entails finding
  special coordinates in which the hypothesis holds.  The conclusion of the 
  main result (Theorem 4.2) is, however, coordinate free. One could  strengthen 
  the   Restriction Hypothesis  so that it is also coordinate free, and this 
  might make a more pleasing statement.  However, it would
  make applications needlessly more difficult.  In most cases the right choice of 
  coordinates is pretty obvious.

\medskip
In Appendix A we present some elementary examples where restriction fails.

In Appendix B certain important algebraic properties of the restriction of quadratic forms are studied.
In particular, Theorem B.9 implies that in geometric cases (where a subset $\GG$ of the 
bundle $G(p,TZ)$ of tangent $p$-planes on a riemannian manifold $Z$ determines the
subequation $F_\GG$), if the submanifold $X$ is totally geodesic, then the restricted subequation
$H\equiv \overline{i^*F_{\GG}}$ on $X$ is geometrically determined by $\GG(TX)$, the
tangential part of $\GG$ along $X$.  That is, $H\equiv \overline{i^*F_{\GG}} = i^*F_{\GG(TX)}$.

In particular, the case $\GG(TX) =\emptyset $ ($X$ is $\GG$-free) is exactly the case when 
$\overline{i^*F_{\GG}} = J^2(X)$, which is uninteresting for restriction   since   
$\overline{i^*F_{\GG}}$ imposes no constraint.  However, this is the appropriate setting for
extension results.

Finally, in Appendix B we give a euclidean example of a subequation
$F_\GG \ss\Sym(\bbr^3)$ and a plane $W\ss \bbr^3$, where $i^*F_\GG$ is not a 
closed set, so that $i^*F_\GG \neq F_{\GG(W)}$.

\medskip
\noindent
{\bf Appendix C. (Extension Theorems).} Intimately related to  restriction is the question of extension, namely,
which functions on a submanifold can be extended to $F$-subharmonic functions in a neighborhood?
In Appendix C we give conditions under which every $C^2$-function has this property.

%\vfill\eject
\vskip .3in

\noindent{\headfont \BB.\  Nonlinear Potential Theory}
\medskip

Suppose $u$ is a real-valued function of class $C^2$ defined on an open subset $X\ss\rn$.
The {\bf full second derivative} or {\bf 2-jet} of $u$ at a point $x\in X$ will be denoted by 
$$
J_x u  \ =\ \left(u(x), D_x u, D^2_x u\right)
\eqno{(\BB.1)}
$$
where   $D_x u = ({\partial u\over \partial x_1}(x),...,{\partial u\over \partial x_n}(x))$
and $D_x^2u  =  (({\partial^2 u\over \partial x_i\partial x_j} (x) ))$.
Occasionally $D^2_xu$ is denoted by $\Hess_xu$.

In this paper constraints on the full second derivative of a function $u\in C^2(X)$ will take the form 
$$
J_x u \ \in\ F_x
\eqno{(\BB.2)}
$$
where $F\ss J^2(X)$ is a subset of the 2-jet space $J^2(X) = X\times \bbr\times \bbr^n\times \Symn$
and $F_x$ denotes the fibre of $F$ at $x\in X$.  Such functions $u$ will be called
{\bf $F$-subharmonic}.

Given an upper semi-continuous 
functions  $u$  on $X$ with values in $[-\infty, \infty)$,   a {\bf test function for $u$ at $x_0$}
is a $C^2$ function $\vf$ defined near $x_0$ which satisfies:
$$
\left.
\cases
{ u-\vf  \ &$\leq$ \ 0 \    \quad {\rm near}\ $x_0$   \cr 
   \ &= \ 0\ \qquad {\rm at}\ $x_0$  
 } 
 \right\}.
\eqno{(\BB.3)}
$$

\Def{\BB.1} An upper semi-continuous  function $u$ on $X$  is {\bf $F$-subharmonic} if for all $x_0\in X$ 
$$
 J_{x_0} \vf \ \in\ F_{x_0} \quad {\rm for\ all\
test\  functions\ } \vf \ {\rm for\ } u \ {\rm at\ } x_0
$$
Let $F(X)$ denote the space of all $F$-subharmonic functions on $X$.
\medskip

Note that if $u(x_0) = -\infty$, then there are no test functions for $u$ at $x_0$.

If $\vf$ is a test function for $u$ at $x_0$, then
so is $\psi\equiv \vf+\half\bra {P(x-x_0)}{x-x_0}$ for any matrix $P\geq0$.
Moreover, $J_{x_0}\psi=J_{x_0}\vf+P$. Consequently, $F(X)$ is empty (except for $u\equiv-\infty$) unless
$F$ satisfies the following {\bf positivity condition (P)} 
$$
F_x + \cp \ \ss\ F_x \fa x\in X
\eqno{(\BB.4)}
$$
where $\cp \equiv \{0\}\times \{0\}\times \{P\in \Symn: P\geq0\}$.
We will abuse notation and also let $\cp$ denote the subset of $\Symn$ of matrices
$P\geq0$.

Assuming this condition (P), it is easy to show that each $C^2$-function $u$ satisfying (\BB.2)
is $F$-subharmonic on $X$.  (The converse is true without (P) since $\vf=u$ is a test function.)

Definition \BB.1 can be recast in a more useful form. (See [HL$_6$,  Lemma 2.4 and Prop. A.1 (IV)].)

\Lemma{\BB.2}  {\sl  Suppose $F\ss J^2(X)$ is a closed  subset, and let $u$ be an upper
semi-continuous function on $X$. Then  $u\notin F(X)$ if and only if 
$\exists\, x_0\in X$, $\a >0$ and $(r,p,A)\notin F_{x_0}$ with}
$$
\eqalign
{
u(x)-\bigl[r+\langle p, x-x_0\rangle +\half \langle A(x-x_0), x-x_0 \rangle\bigr] \ &\leq \ -\a|x-x_0|^2 \qquad
{\rm near\ } x_0 \qquad{\sl and}  \cr
&=\ 0\qquad\qquad\qquad\ \ \ {\sl at}\ \ x_0
}
$$

Using this Lemma, basic potential theory for $F$-subharmonic functions is elementary to establish.
See  Appendices A and B in [HL$_6$].

\Theorem{\BB.3} {\sl  
Let $F$ be an arbitrary  closed subset of $\Jt(X)$.
\medskip

\item{(A)} (Local Property)  $u$ is locally  $F$-subharmonic if and only if $u$ is globally
 $F$-subharmonic.

\medskip

\item{(B)}  (Maximum Property)  If $u,v \in F(X)$, then $w=\max\{u,v\}\in F(X)$.

\medskip

\item{(C)}     (Coherence Property) If $u \in F(X)$ is twice differentiable at $x\in X$, then $\jtx u\in F_x$.

\medskip

\item{(D)}  (Decreasing Sequence Property)  If $\{ u_j \}$ is a 
decreasing ($u_j\geq u_{j+1}$) sequence of \ \ functions with all $u_j \in F(X)$,
then the limit $u=\lim_{j\to\infty}u_j \in F(X)$.

\medskip

\item{(E)}  (Uniform Limit Property) Suppose  $\{ u_j \} \ss F(X)$ is a 
sequence which converges to $u$  uniformly on compact subsets to $X$, then $u \in F(X)$.

\medskip

\item{(F)}  (Families Locally Bounded Above)  Suppose $\cf\subset F(X)$ is a family of 
functions which are locally uniformly bounded above.  Then the upper semicontinuous
regularization $u=v^*$ of the upper envelope 
$$
v(x)\ =\ \sup_{f\in \cf} f(x)
$$
belongs to $F(X)$.
}
\medskip

There are certain obvious additional properties (e.g. {\sl If $F_1\ss F_2$, then $u\in F_1(X)\ \Rightarrow\ 
u\in F_2(X)$}), which will be used without reference.

Although the positivity condition (P) is not needed in the proofs of either Lemma \BB.2 or Theorem \BB.3,
the fact that without (P) there are no $F$-subharmonic functions, other than $u\equiv-\infty$,
explains this requirement.

\Def{\BB.4}  A closed subset $F\ss J^2(X)$ which satisfies the positivity condition (P)
will be called a {\bf subequation}.

\medskip
\noindent
{\bf Note.} This does not agree with the terminology of [HL$_6$] where subequations were
assumed to have two additional properties: a stronger topological condition (T) 
and, in order to have a chance of proving uniqueness in the Dirichlet problem, 
 standard negativity condition (N) on the values of the dependent variable (cf. [HL$_6$]).
However, these conditions are unnecessary for the discussion in this paper.
\medskip

The following basic example will be elaborated later in Examples \FF.2 and \MM.7.

\Ex{\BB.5.  (The Monge-Amp\`ere equation $\det(D^2 u)=0$)}  There are $n$ 
different subequations (or branches) associated with this equation.  
Thus it generates $n$ distinct notions
of subharmonic.  The $q^{\rm th}$ branch, denoted here by $\L_q$,
 is defined by the inequality $\l_q\geq0$,
where $\l_{\rm min} (A)= \l_1(A) \leq \cdots \leq \l_n(A)=\l_{\rm max}(A)$
are the ordered eigenvalues of $A\in\Symn$.
Equivalently, $D_x^2u$ (or $D_x^2\vf$ with $\vf$ a test function for $u$ at $x$)
is required to have at least $n-q+1$ eigenvalues which are $\geq0$.
Note that $\L_{\rm min} = \bbr\times\rn\times\cp$  is the smallest branch.  It
follows easily from Lemma \BB.2 that classical convex functions
are $\L_{\rm min}$-subharmonic.  The converse is also true, 
but the proof does require the restriction theorem and may be 
 considered its  most elementary application (see Example \CC.3).

The largest branch $\L_{\rm max}$, where only one eigenvalue is required to be $\geq0$
is particularly important.  The $\L_{\rm max}$-subharmonic functions can be described more
concretely using a class of functions introduced in  [HL$_4$].

\Def{\BB.6} A function $u\in \USC(X)$ is said to be {\sl subaffine on $X$}
if for each compact subset $K\ss X$ and each affine function $a$,
$$
u\ \leq\ a \ \ \ {\rm on}\ \ \partial K 
\qquad\Rightarrow\qquad
u\ \leq\ a \ \ \ {\rm on}\ \ K 
\eqno{(\BB.5)}
$$

In   [HL$_4$, Remark 4.9] we proved the following.

\Prop{\BB.7} {\sl
Given  $u\in \USC(X)$, the following are equivalent.

\medskip
(1)\ \ \ $u$\ is locally subaffine,
\medskip
(2)\ \ \ $u$\ is $\L_{\rm max}$-subharmonic on $X$,
\medskip
(3)\ \ \ $u$\ is  subaffine on $X$.
}

\medskip

For the sake of completeness we give a different, shorter proof here.
\smallskip
\noindent
{\bf Proof that (1) \ $\Rightarrow$ \ (2):}
Suppose $u$ is not $\L_{\rm max}$-subharmonic on $X$.
Apply Lemma \BB.2.
Since $\l_{\rm max}(A) \geq 0$ is false if and only if $A<0$, it follows directly
 that $u$ is not sub-the-affine-function $r+\bra p {x-x_0} $
on small balls about $x_0$.

\smallskip
\noindent
{\bf Proof that (2) \ $\Rightarrow$ \ (3):}
Suppose $u$ is not subaffine on $X$.
Then there exists a compact set $K\ss X$ and an affine function $a$
such that (\BB.5) fails, that is, $u-a$ has a strict interior maximum on $K$.
Thus,  for $\e>0$ sufficiently small, the function $u(x) +{\e\over 2}|x|^2-a(x)$ 
also attains its maximum value (say $k$) at an interior point $x_0$ of $K$. 
Now the function $\vf(x) = a(x) -{\e\over 2}|x|^2 +k$ is a test function for 
$u$ at $x_0$. Since $D_{x_0}^2\vf = -\e I$, $u$ is not $\L_{\rm max}$-subharmonic
0n $X$.\qed

%%%%%%%%%%%%%%%%%%%%%%%%%%%%%%%%%%%%%%%%%%%%%%%%
%%%%%%%%%%%%%%%%%%%%%%%%%%%%%%%%%%%%%%%%%%%%%%%%
%%%%%%%%%%%%%%%%%%%%%%%%%%%%%%%%%%%%%%%%%%%%%%%%
%%%%%%%%%%%%%%%%%%%%%%%%%%%%%%%%%%%%%%%%%%%%%%%%
%%%%%%%%%%%%%%%%%%%%%%%%%%%%%%%%%%%%%%%%%%%%%%%%

%\vfill\eject
\vskip.3in

\noindent{\headfont \CC. An Introduction to Restriction -- The Geometric Case in $\rn$.}

\bigskip

In this section we describe a special case of our restriction results which is simple but important.
A subequation $F$ is said to be {\bf geometrically determined} by a closed subset 
$\GG$ of the Grassmannian  $G(p,\rn)$ of (unoriented) $p$-planes through the origin in $\rn$ if 
$F\equiv F_\GG$ is defined by
$$
\trace \left\{ D^2_x u\bigr|_W\right\}\ \geq\ 0 \fa W\in \GG
\eqno{(\CC.1)}
$$
and for all $x\in X$.  The upper semi-continuous functions in $F_\GG(X)$ will be
referred to as $\GG$-plurisubharmonic on $X$.

\Ex{\CC.1. (Classical Subharmonicity)}
If $p=n$ and $\GG=G(n,\rn)=\{\rn\}$, then $u$ is $\GG$-plurisubharmonic on the open set
$X\ss\rn$ if and only if $u$ is subharmonic ($\trace ( D^2 u) = \Delta u\geq0$ in the $C^2$-case) using
any of the equivalent classical definitions ($u\equiv -\infty$ on components of $X$ is allowed).
In the case $n=1$, subharmonicity is the same as classical convexity in one variable,
expanded to allow $u\equiv -\infty$  as a matter of convenience.
\medskip

An {\bf affine $\GG$-plane} is an affine plane in $\rn$ whose translate through the origin belongs to $\GG$.

\medskip
\noindent
{\bf Restriction Theorem \CC.2.} {\sl  A function $u$ is $\GG$-plurisubharmonic on $U\ss\rn$ if and only if}
$$
u\bigr|_{U\cap W} \ \ {\rm is\ subharmonic\ for\ each\ affine\ } \GG \, {\rm plane\ } W.
\eqno{(\CC.2)}
$$

\pf
Half of the proof is trivial.  If $\vf$ is a test function for $u$ at $x_0\in X$ with $J_{x_0} \vf \notin F_{x_0}$,
then by definition of $F\equiv F_\GG$ there exists a $W\in\GG$ with $\tr_W D^2_{x_0} \vf <0$.
Therefore (cf. Ex. \CC.1) $u\bigr|_{X\cap (W+x_0)}$ is not subharmonic at $x_0$.
The other half, namely the assertion  that restrictions of $\GG$-psh functions to affine $\GG$-planes are subharmonic
is proved in the   Section \FF.  It is a special case of  our general Geometric Restriction Theorem
\MM.2.\qed

\Ex{\CC.3. (Classical Convexity)}  If $\GG=G(1,\rn)$, then this restriction theorem is precisely the theorem
required to establish that  the condition $D^2u\geq 0$ in the viscosity sense implies that  $u$ is convex
(or possibly $\equiv -\infty$).  Somewhat surprisingly we were unable to find an elementary viscosity proof
of this fact in the literature.  Such  a proof is essentially given in [HL$_4$, Prop. 2.6], and this is the prototype
of our proof of the general restriction theorem.

\Ex{\CC.4. (Plurisubharmonicity in Complex Analysis)}  A function $u\in \USC(X)$ with $X$ an open subset of
$\bbc^n$ is said to be plurisubharmonic if the restriction of $u$ to each affine complex line is classically subhharmonic. Our Restriction Theorem \CC.2 states that this classical notion is equivalent to being $\GG$-plurisubharmonic where $\GG=G_\bbc(1, \bbc^n) \ss G_\bbr(2,\bbc^n)$ is the Grassmannian of complex lines in $\bbc^n$.
\medskip

Further examples abound.  A wide class (including Examples \CC.3 and \CC.4)
is given by choosing a calibration $\phi\in \L^p\rn$ and then setting 
$$
\GG(\phi) \ \equiv\ \left\{  W\in G(p,\rn):  \phi\bigr|_W   \ {\rm is\ the\ standard\ volume\ form\ on\ } W   \right\}
\eqno{(\CC.3)}
$$
for one of the choices of orientation on $W$.

%%%%%%%%%%%%%%%%%%%%%%%%%%%%%%%%%%%%%%%%%%%%%%%%
%%%%%%%%%%%%%%%%%%%%%%%%%%%%%%%%%%%%%%%%%%%%%%%%
%%%%%%%%%%%%%%%%%%%%%%%%%%%%%%%%%%%%%%%%%%%%%%%%
%%%%%%%%%%%%%%%%%%%%%%%%%%%%%%%%%%%%%%%%%%%%%%%%
%%%%%%%%%%%%%%%%%%%%%%%%%%%%%%%%%%%%%%%%%%%%%%%%

\vskip.3in
%\vfill\eject

\noindent{\headfont \EE.   The General  Restriction Theorem.}

\bigskip

Suppose $Z$ is an open subset of $\bbr^N = \rn\times \bbr^m$ 
with coordinates $z=(x,y)$.   Set $X=\{x\in \rn : (x,y_0)\in Z\}$ for a fixed $y_0$,
and let $i: X \hookrightarrow Z$ denote the inclusion map $i(x) = (x, y_0)$. 
Adopt the notation
\smallskip
\centerline
{ $r= \vf(x,y_0)$, $p={\partial \vf\over \partial x} (x,y_0)$, 
 $q={\partial \vf\over \partial y} (x,y_0)$, 
 $A={\partial^2 \vf\over \partial x^2} (x,y_0)$, 
 $B={\partial^2 \vf\over \partial y^2} (x,y_0)$, 
 $C={\partial^2 \vf\over \partial x\partial y} (x,y_0)$
 }
 \smallskip
 \noindent
 for the 2-jet $J_z\vf$ of a function $\vf$ at $z=(x,y_0)$.
Then the 2-jet of th restricted function $\psi(x) = \vf(x,y_0)$ is given by $J_x\psi = (r,p,A)$.  
Thus, restriction $i^*:J^2(Z)\arr J^2(X)$ on 2-jets is given by
$$
i^*\left(r, \ (p,q), \ \left(\matrix{ A&C\cr C^t&B}\right)\right)\ \ =\ \ (r,p,A) \qquad {\rm at\ \ } i(x)=z.
\eqno{(\EE.1)}
$$ 

If $F$ is a subset of $J^2(Z)$, then  the {\bf restriction $i^*_X F$  of $F$ to $X$} is a subset of
$J^2(X)$.
Each quadratic form $P\geq0$ on $\rn$ is the restriction of a quadratic form $\wt P\geq0$ on $\bbr^N$.
This proves that:
$$
{\rm  If \ } F \ {\rm satisfies\  the\ positivity\  condition\ (P),\ then\ } i^*F \ {\rm \ also\  satisfies\  (P).}
\eqno{(\EE.2)}
$$
We shall also consider the  closure $H =  \overline{i^*F}$.  It is obvious that 
$$
F\ \ {\rm satisfies \ \ (P)}\quad\Rightarrow \quad   H \ \ {\rm satisfies \ \ (P)}
\eqno{(\EE.3)}
$$
Thus, $H\equiv \overline{i^*F}$ is a subequation (Def. \BB.4), and it will be referred to
as the {\sl restricted subequation}.

\Def{\EE.1} We say that {\bf restriction to $X$ holds for $F$} if 
$$
u \  {\rm  is\ } F  {\rm -subharmonic \ on \ } Z\quad\Rightarrow \quad 
u \ \bigr|_X   \ {\rm  is\ } H {\rm -subharmonic \ on \ } X
\eqno{(\EE.4)}
$$

This is not always the case. 
Some elementary  examples are presented in Appendix A.
 Of course, if $u\in C^2(Z)$ is $F$-subharmonic, then $u\bigr|_X$ is $H$-subharmonic 
 on $X$ since $i^*Ju = J i^*u$. The only issue is with $u\in\USC(Z)$ that are not $C^2$.
 Let $\Jtn = J^2_0(\rn) = \bbr\oplus \rn\oplus\Symn$.

\medskip
\noindent
{\bf The Restriction Hypothesis:}  
Given $x_0\in X$ and  $(r_0, p_0, A_0)\in \Jtn$ and given
$z_\e = (x_\e, \, y_\e)$ and $r_\e$ for a sequence of real numbers $\e$ converging to 0. \medskip
$$
{\rm If\ \ } \left( r_\e,\ \left(p_0+A_0(x_\e-x_0), \ {{y_\e-y_0}\over\e}\right), \ \left(\matrix{A_0&0\cr 0&{1\over \e}I}\right)\right)\ \ \in\ \ F_{z_\e} 
\eqno{(\EE.5)}
$$
$$
{\rm and\ \ }  x_\e\ \to\  x_0,\ \ {{|y_\e-y_0|^2}\over\e}\ \to\ 0, \ \  r_\e\ \to\ r_0,
\eqno{(\EE.6)}
$$
then
$$
(r_0, p_0, A_0)\ \in \ H_{x_0}.
$$

\medskip
\noindent
{\bf Remark.}  If the subequation $F$ is independent of the $r$-variable, that is, if 
$F_x$ can be considered as a subset of the reduced 2-jet space $\overline{\bbj}^2_{x} =\rn\times \Symn$,
then the restriction hypothesis can be restated as follows.

\medskip
\noindent
{\bf Restriction Hypothesis (Second Version -- for $r$-Independent  Subequations):}  Given $x_0\in X$ and $z_\e = (x_\e, y_\e)$ converging to
$z_0=(x_0,y_0)$ with ${1\over \e}|y_\e-y_0|^2 \to0$, for a sequence of real numbers $\e$ converging to 0, 
consider the polynomials
$$
\psi_\e(x,y)\ \equiv\    r_0+\langle p_0, x-x_0\rangle +\half \langle A_0(x-x_0), x-x_0 \rangle + {1\over 2\e}|y-y_0|^2.
\eqno{(\EE.7)}
$$ 
$$
{\rm If\ \ } \overline{J}_{z_\e}\psi_\e  \in  F_{z_\e} \fa \e, \ \ {\rm then\ \ } (p_0,A_0) \in H_{z_0}
$$
\medskip

This follows since the reduced jet $ \overline{J}_{z_\e}\psi_\e$ equals the jet in (\EE.5) modulo $r_\e-r_0$.

\smallskip
\noindent
{\bf The General Restriction Theorem \EE.2.}  {\sl  Suppose $u\in \USC(Z)$.  Assume the restriction hypothesis. Then with $H\equiv \overline{ i^*F}$,}
 $$
 u\in F(Z)  \qquad \Rightarrow\qquad u\bigr|_X \in H(X).
 $$

\Remark{\EE.3}  See Example B.6 in Appendix B for a case where $i^*F$ is not closed.

\pf
If $u\bigr|_X \notin H(X)$, then by Lemma \BB.2 (since $H$ is closed) there exists $x_0\in X$, $\a>0$, and $(r_0, p_0, A_0)\notin H_{x_0}$ such that 
$$
\eqalign
{
u(x,y_0) - Q(x)\ &\leq \ -\a|x-x_0|^2 \qquad
{\rm near\ } x_0 \qquad{\rm and}  \cr
&=\ 0\qquad\qquad\qquad\ \ \ {\rm at}\ \ x_0
}
\eqno{(\EE.8)}
$$
where 
$$
 Q(x) \ \equiv\  r_0+\langle p_0, x-x_0\rangle +\half \langle A_0(x-x_0), x-x_0 \rangle.
\eqno{(\EE.9)}
$$

In the next step we construct $z_\e=(x_\e, y_\e)$ satisfying (\EE.6) with $r_\e \equiv u(z_\e)$.
Set
$$
w(x,y)\ \equiv\ u(x,y) - Q(x).
$$
Let  $B(z_0)$ denote a small closed ball about $z_0$ in $\bbr^N$, so that  (\EE.8) holds
on the $y_0$-slice.
For each $\e>0$ small, let
$$
M_\e\ \equiv \ \sup_{B(z_0)} \left(w - \smfrac1{2\e} |y-y_0|^2\right),
\eqno{(\EE.10)}
$$
and choose  $z_\e$ to be a maximum point.  Since the value of this function at $z_0$ is zero,
the maximum value $M_\e\geq0$.  Furthermore, the $M_\e$ decrease to a limit, say $M_0$.
Now
$$
\eqalign
{
M_\e\ &=\ w(z_\e) - \smfrac1{2\e}|y_\e-y_0|^2\ =\ 
w(z_\e) - \smfrac1{4\e}|y_\e-y_0|^2- \smfrac1{4\e}|y_\e-y_0|^2  \cr
&\leq M_{2\e}- \smfrac1{4\e}|y_\e-y_0|^2,\qquad {\rm that\ is}   \cr
}
$$
$$
 \smfrac1{4\e}|y_\e-y_0|^2  \ \leq \  M_{2\e}- M_{\e}.
$$
Thus
$$
\smfrac1{\e}|y_\e-y_0|^2 \ \arr\ 0
\eqno{(\EE.11)}
$$
and in particular $y_\e\to y_0$.

Suppose now that $\bar z=(\bar x, y_0) $ is a cluster point of $\{z_\e\}$.  Then taking a sequence
$z_\e \to \bar z$
$$
M_0\ =\ \lim_{\e\to0} M_\e \ =\ \lim_{\e\to0}(w(z_\e) - \smfrac1{2\e}|y_\e-y_0|^2)
\ =\ \lim_{\e\to0}w(z_\e)\ \leq\ w(\bar z)
\eqno{(\EE.12)}
$$
by (\EE.10), (\EE.11) and the fact that $w$ is upper semi-continuous.  By (\EE.8) and the fact that $\bar y=y_0$,
we have  $w(\bar z)\leq0$.   Hence, $M_0=w(\bar z)=0$.
Since $w(x,y_0)$ has a strict maximum of 0 at $z_0=(x_0,y_0)$, and this maximum value is 
attained at $\bar z=(\bar x, y_0)$, we must have $\bar x=x_0$.  Thus
$$
x_\e\ \to\ x_0.   % \qquad {\rm and}\qquad z_\e\ \in\ 
%\Int B(z_0)\quad{\rm for \ \ } \e>0\ \ {\rm small}.
\eqno{(\EE.13)}
$$
Now by (\EE.12), we have $0 = \lim_{\e\to 0} w(z_\e) = 
\lim_{\e\to 0}  \bigl(u(z_\e)-Q(z_\e) \bigr)   = \lim_{\e\to 0} r_\e -r_0$,
which completes the proof that (\EE.6) is satisfied.

It remains to verify (\EE.5).  The notation has been arranged so that 
$$
u-\psi_\e\ =\ w-{1\over 2\e}|y-y_0|^2
\eqno{(\EE.14)}
$$
where $\psi_\e$ is defined by (\EE.7).  Consequently, (\EE.10) can be restated as 
$$
\eqalign
{
u-\psi_\e  \ &\leq \ M_\e \qquad  {\rm near\ }  z_\e \qquad{\rm and}  \cr
&=\ M_\e \qquad {\rm at}\ \ z_\e,
}
\eqno{(\EE.10)'}
$$
that is, $\vf_\e \equiv \psi_\e +M_\e$ is a test function for $u$ at $z_\e$.
This implies that $\jt_{z_\e}\vf_{\e} \in F_{z_\e}$.  Computing this 2-jet verifies (\EE.5).
The Restriction Hypothesis now  implies that $(r_0,p_0,A_0)\in H_{x_0}$,
which is a contradiction.\qed

%%%%%%%%%%%%%%%%%%%%%%%%%%%%%%%%%%%%%%%%%%%%%%%%
%%%%%%%%%%%%%%%%%%%%%%%%%%%%%%%%%%%%%%%%%%%%%%%%
%%%%%%%%%%%%%%%%%%%%%%%%%%%%%%%%%%%%%%%%%%%%%%%%
%%%%%%%%%%%%%%%%%%%%%%%%%%%%%%%%%%%%%%%%%%%%%%%%
%%%%%%%%%%%%%%%%%%%%%%%%%%%%%%%%%%%%%%%%%%%%%%%%

\vfill\eject

\noindent{\headfont \FF.   First Applications.}

\medskip

We now examine some applications of the Restriction Theorem \EE.2
\bigskip
\centerline{\bf Restriction in the Constant Coefficient Case}
\medskip

Suppose $F=Z\times {\bf F}$ for ${\bf F}\ss  \JtN$.  Then $F$ is said to have 
{\bf constant coefficients} on $Z$.  Now consider $X=  Z \cap \{ y=y_0\}$ as above.
If $F$ has constant coefficients on $Z$, then the restriction of 2-jets 
gives a set $H=\overline{ i^*F} = X\times {\bf H}$ with constant coefficients on $X$.

\Theorem {\FF.1. (Restriction for Euclidean Subequations)}
{\sl  Suppose $F\ss J^2(Z)$ is closed,  has constant coefficients and satisfies (P).  Then
\medskip
\centerline
{
$u$ is $F$-subharmonic on $Z$ \qquad$\Rightarrow\qquad u\bigr|_X$ is $H$-subharmonic on $X$.
}
}

\pf
In this case the restriction hypothesis is easy to verify.  Since 
$$
\left( r_\e, (p_\e, q_\e), \left(\matrix{A_0&0\cr 0&\smfrac 1 \e I}\right)\right) \ \in\ F_{z_e} \ =\ {\bf F},
$$
we have that  the restricted 2-jet $( r_\e, p_\e, A_0) \in H$ even though $z_\e \notin X$.
Now the fact that  $r_\e\to r_0$ and $p_\e = p_0 +A_0(x_\e-x_0) \to p_0$ is enough to conclude
 that $(r_0,p_0,A_0)\in  {\bf H}=   H_{z_0}$.\qed

\medskip

There are many subequations  for which Theorem \FF.1  is interesting. 
For one such basic case we continue with Example \BB.5.

\Ex{\FF.2. (Branches of the Homogeneous Monge-Amp\`ere Equation)} 
The $q^{\rm th}$ branch $\L_q$  of the homogeneous Monge-Amp\`ere  equation
on $\rn$ is defined by requiring that the $q^{\rm th}$ ordered eigenvalue of the second derivative
be $\geq0$, i.e., the subequation $\L_q$  is defined by
$$
\l_q(A) \ \geq\ 0 \qquad{\rm for}\ \ A\in\Symn.
\eqno{(\FF.1)}
$$
Even though this is not one of the geometric cases, the 
subharmonics can be characterized via restriction,
providing an extension of Proposition \BB.7.

\Theorem{\FF.3} {\sl
A function $u\in \USC(X)$ is $\L_q$-subharmonic if and only if its 
restriction to each affine $q$-plane $V\ss\rn$ is subaffine (see Definition  \BB.6).
}

\pf
In order to apply the Restriction Theorem \FF.1 to a $\L_q$-subharmonic function
on $\rn$ we must first compute the restricted subequation on an affine $q$-plane $V$.
We can assume that $V$ is a vector subspace of $\rn$.  Given 
$A\in\Symn$, recall that:
$$
\l_q(A)\ =\ \inf_{W} \l_{\rm max}\left(  A\bigr|_W\right)
\eqno{(\FF.2)}
$$
where the inf is taken over all $q$-dimensional subspaces $W\ss \rn$,
and  $ \l_{\rm max}\left(  A\bigr|_W\right) =  \l_{q}\left(  A\bigr|_W\right)$.
It follows that the subequation $\L_q$ on $\rn$ restricts to the 
 subequation $\L_q$ on $\bbr^p$ for any $p\geq q$.  Now on $\bbr^q$,
 $  \l_{q}(B) = \l_{\rm max}(B)$ so that $\l_q(B)\geq0$ on $\bbr^q$
 if and only if at least one eigenvalue of $B$ is $\geq 0$.  
 Combining the  Restriction Theorem \FF.1  with Proposition \BB.7 completes the 
 proof in one direction.

If $u$ is not $\L_q$-subharmonic on $\rn$, then using  Lemma \BB.2 and  some
normalizations, one sees that  there exists $A$ with $\l_q(A) <0$ 
such that $u(x) - \bra {Ax}x \leq 0$ near $x=0$ with equality at $x=0$.
Take  $V$ to be the span of the first $q$ ordered eigenvectors of $A$.
Then $u\bigr|_V -A\bigr|_V \leq 0$ near $x=0$  and $A\bigr|_V<0$, proving that
$u\bigr|_V $ is not subaffine.   \qed

\Remark{\FF.4}  This theorem easily extends to subequations defined
by $\l_q(A) \geq f(r, |p|)$ with $f(r,s)$ non-decreasing in $s$ and continuous in $r$.

\Ex{\FF.5. ($i^*F$ not Closed)}
Define $F$ on $\bbr^2$ by $|p||q|\geq1$.  Then $i^*F$  on $\{y=0\}$ is defined 
by $p\neq0$, and $H=\overline{ i^*F}$ is all of $J^2(\bbr)$.  In particular, $i^*F$ is not closed.  
A more interesting (geometrically defined) example where $i^*F$ is not closed, is given in
Appendix B.

%%%%%%%%%%%%%%%%%%%%%%%%%%%%%%%%%%%%%%%%%%

%\vfill\eject 
\bigskip

\centerline{\bf The Geometric Case in $\rn$.}
\medskip

As in Section \CC \ suppose that $F_\GG$ is geometrically defined by  closed subset
$\GG$ of the grassmannian $G(p,\bbr^N)$.

\medskip
\noindent
{\bf Proof of Theorem \CC.2.}  It  is a special case of Theorem \FF.1.
To see this suppose ${\bf W}$ is an affine $\GG$-plane with (constant) tangent plane $W\in\GG$.
Then for any quadratic form $Q$ at any point of ${\bf W}$ we have $\tr_W i^*_{\bf W} Q = \tr_W Q$
which proves that $ i^*_{\bf W} \bbf_\GG \ss \bbf_{\{W\}}$,  the classical 
(subharmonic) subequation on ${\bf W}$ (cf. Example \CC.1). \qed\medskip

This Restriction Theorem \CC.2 can be 
generalized by considering a subspace $V\ss\bbr^N$ of larger dimension $n\geq p$
and defining 
$$
\GG(V) \ \equiv\ \{W\in \GG : W\ss V\}
\eqno{(\FF.3)}
$$
to be the space of $\GG$-planes which are {\bf tangential to $V$}.  Since $\GG(V)$ is a 
closed subset of the grassmannian $G(p, \bbr^N)$, it geometrically determines a subequation
$F_{\GG(V)}$ on $V$ by  
$$
F_{\GG(V)}\ \equiv\ \{a\in \Sym(V^*) : \tr_W a\geq 0 \ \ \forall \ W\in \GG(V)\}
\eqno{(\FF.4)}
$$

\Theorem{\FF.6} {\sl
If $u$ is $\GG$-plurisubharmonic on an open subset $U\ss \bbr^N$, then
for each affine subspace $V$ of $\bbr^N$,
$$
u\bigr|_{U\cap V} \ \ {\rm is}\ \ \GG(V) \ {\rm plurisubharmonic}.
$$
}
\medskip

  Theorem \CC.2 is the special case where   $V=W$ and so   $\GG(V) =\{W\}$.

\Remark{\FF.7}  As in Theorem \CC.2 the converse (where one considers all affine subspaces
$V$ of dimension $n$ with $n\geq p$) is trivial.

\pf 
Let $i^*_V$ denote the restriction of 2-jets from $\bbr^N$ to $V=\rn$.  For $W\ss V$ one has 
$\tr_W i^*_VQ = \tr_W Q$ for all quadratic forms $Q$, which proves that
$$
i^*_V F_\GG \ \ss \ F_{\GG(V)}.
\eqno{(\FF.5)}
$$
Therefore $\overline {i^*_V F_\GG } \ss  F_{\GG(V)}$,  and so 
Theorem \FF.6 is a special case of Theorem \FF.1.\qed   
\medskip

  In Appendix B  (Theorem B.3) we prove that in fact $F_{\GG(V)}$ is the restricted subequation,
  i.e., 
  $$
\overline{  i^*_V F_\GG } \ = \  F_{\GG(V)}.
  $$

%Theorem \FF.1 combined with Proposition \FF.2 
%yields a generalization of Theorem \CC.2.

%%%%%%%%%%%%%%%%%%%%%%%%%%%%%%%%%%%%%%%%%%

\vfill\eject
%\bigskip

\centerline{\bf Subequations which can be 
Defined Using Fewer of the Variables in $\bbr^N$.}
\medskip

Suppose that $F$ can be defined using fewer of the variables in $\bbr^N$, say
using only the variables in $\rn\ss\bbr^N$.  This means by definition that  there exists $\bbh \ss \Jtn$
with $\bbf = (i^*)^{-1} \bbh$ where 
$i^*:\JtN\to \Jtn$ is the restriction map.

We shall say that a  function $u\in \USC(Z)$ is {\sl horizontally $H$-subharmonic on an open set $Z\ss \bbr^N$}
if for each $y_0\in \bbr^n$ the function $u(x,y_0)$ is of type $H$ on $Z\cap \{y=y_0\}$.
\medskip

As another  special case of Theorem \FF.1 we have

\Theorem{\FF.8}  {\sl  Suppose the constant coefficient subequation 
$\bbf = (i^*)^{-1}(\bbh)$ can be defined using   
 the variables $\rn\ss\bbr^N$.  Then $u$ is $F$-subharmonic on $Z$ if and only if 
 $u$ is horizontally $H$-subharmonic on $Z$.}

\bigskip

\centerline{\bf Families of Subequations.}
\medskip

Theorem \FF.8 extends to  a more general, non constant coefficient situation.
Let $F(y) \ss J^2(X)$ be a family of subequations parameterized by 
points $y$ in an open subset $Y\ss \bbr^m$.  Consider the subset
$F\ss J^2(Z)$, $Z\equiv X\times Y$, defined by 
$$
J^2_z \vf (z) \ \in\ F_z\qquad\iff\qquad J^2_x \vf(x,y) \ \in \ F_x(y)  \quad z=(x,y)
\eqno{(\FF.5)'}
$$
Obviously, $F$ satisfies the positivity condition (P).  Note that $F$ is a subequation 
in the sense of Definition \BB.4 if and only if
$F\ss J^2(Z)$ is closed. In this case we say the family  $\{F(y)\}$ is  {\sl closed}.

\Theorem{\FF.9}  {\sl  Suppose $\{F(y)\}$ is a closed family of subequations as above.
Then a function $u \in\USC(Z)$ is $F$-subharmonic if and only if the restriction $u(x,y_0)$ is
$F(y_0)$-subharmonic on $X$ for each $y_0\in Y$.
}

\pf
If $\vf(x,y)$ is a test function for $u(x,y)$ at $z_0=(x_0,y_0)$, then $\vf(x,y_0)$ is a test function
for $u(x,y_0)$ at $x_0$. If $u(x,y_0)$ is $F(y_0)$-subharmonic, then 
$J^2_{x_0} \vf(x_0,y_0) \in F_{x_0}(y_0)$, or equivalently, $J^2_{z_0}\vf \in F$.

Conversely, assume $u$ is $F$-subharmonic on $Z$.  Consider the data in the restriction
hypothesis.  By the definition (\FF.5)$'$ of $F$, the condition (\DD.5) can be restated as 
$$
J_\e \ =\ (r_\e, p_0+A_0(x_\e-x_0), A_0)\ \in\ F_{x_\e}(y_\e).
$$
Since $z_\e\to z_0$ and $F$ is closed, this implies that $(r_0, p_0, A_0) = \lim J_\e$
must belong to $F_{x_0}(y_0)$. The result now follows from Theorem \EE.2.\qed

%%%%%%%%%%%%%%%%%%%%%%%%%%%%%%%%%%%%%%%%%%%%%

\bigskip
%\vfill\eject

\def\bll{I \!\! L}

\centerline{\bf Restriction in the Linear Case}
\medskip

Consider the second-order linear operator with smooth coefficients:
$$
\bll\left(z, r, (p, q), \left(\matrix{A&C\cr C& B}\right)\right)  \ \equiv\ 
\langle a(z), A \rangle + \langle \a(z), p \rangle +\g(z) r + 
 \langle b(z), B \rangle + \langle \b(z), q \rangle + \langle c(z), C \rangle  
$$
Let $L\ss Z\times \JtN$ be the subset defined by $\bll\geq0$.
Then, of course, $L$ is a subequation (i.e., positivity holds) if and only if
$$
\left(
\matrix
{
a(z) & c(z) \cr
c(z) & b(z)
}
\right)
\ \geq\ 0,
$$
in which case $L$ will be referred to as a {\sl linear subequation}.
Consider
$H_{x}\equiv i^*L_{z}$ with $z=(x,y_0) \in X$.

We will prove that restriction holds in two cases, which taken together ``essentially'' exhaust the
linear operators $\bll$.  In the first case we assume that at least on of the coefficients $\b(x_0, y_0), b(x_0, y_0)$
or $c(x_0, y_0)$ in non-zero.  Restriction locally holds but is completely trivial since
$H_x = \bbj^2_n$ is everything for $x$ near $x_0$.  If, for example, $\b(x_0, y_0)\neq0$, then
by choosing $q$ to be a sufficiently large multiple of $\b(x_0, y_0)$, any jet $(r,p, A)$
can be shown to lie in $H_x$.

The second case is much more interesting.  We assume the following
{\bf linear restriction hypothesis}:
$$
\b(x, y_0), b(x, y_0),\ {\rm and}\ c(x, y_0)\ \ {\rm vanish\ identically\ on\ } X
\eqno{(\FF.6)}
$$
Define the linear operator
$$
\bll_X(x,r,p A) \ \equiv \bra {a(x, y_0)}{A} + \bra {a(x, y_0)}{p} + \g(x, y_0)r
\eqno{(\FF.7)}
$$
on $X$.  Under this hypothesis $H \equiv i^*F$ is the subset of $X \times\bbj^2_n$ defined by 
the linear inequality $\bll_X\geq0$.

\Theorem{\FF.10}  {\sl  Assume that $L$ is a linear subequation satisfying  the linear restriction hypothesis.  Then
\medskip
\centerline
{
$u$ is $\bll$-subharmonic on $Z$ \qquad$\Rightarrow\qquad u\bigr|_X$ is $\bll_X$-subharmonic on $X$.
}
}

\pf
Since $\b$ vanishes on $X$, we have $|\b(x,y)| \leq C|y-y_0|$.  Moreover, since
$b$ vanishes on $X$ and since (P) implies $b(z)\geq0$, $b$ must vanish to second order, i.e., 
$|b(x,y)|\leq C|y-y_0|^2$.  These two facts are enough to verify the 
restriction hypothesis in Lemma  \EE.1.  Assume that
$$
\eqalign
{
0\   &\leq\ \bll\left(z_\e, r_\e, (p_0+ A_0(x_\e-x_0),\, \smfrac{y_\e-y_0}\e), 
\left(\matrix{A_0&0\cr 0& \smfrac 1 \e I}\right)\right)  \cr
&=\ 
\langle a(z_\e), A_0 \rangle + \langle \a(z_\e), p_0 + A_0(x_\e-x_0)\rangle +\g(z_\e) r_\e + 
 \langle b(z_\e), \smfrac 1 \e I \rangle + \langle \b(z_\e), \smfrac{y_\e-y_0}\e \rangle \cr
 }
$$
and that 
$$
x_\e \ \to \ x_0, \ \ {|y_\e-y_0|^2\over \e}\ \to\ 0, \ \ {\rm and}\ \ r_\e\ \to\ r_0.
$$
Now 
$$
\left|\left\langle \b(z_\e),  {y_\e-y_0\over \e}  \right\rangle \right|
\ \leq \ C\,  {|y_\e-y_0|^2\over \e} \ \to\ 0\quad {\rm and}
$$
$$
\left|\left\langle b(z_\e),  {1\over \e} I \right\rangle \right|
\ \leq \ C\,  {|y_\e-y_0|^2\over \e} \ \to\ 0.
$$
Hence the RHS converges to 
$$
\langle a(z_0), A_0 \rangle + \langle \a(z_0), p_0 \rangle +\g(z_0) r_0
\ =\ \bll_X(z_0,r_0, p_0 A_0) 
$$
which proves that $(z_0,r_0, p_0, A_0) \in H_{x_0}$.\qed

\vfill\eject

\Remark{\FF.11. (Versions of the Linear Restriction Hypothesis)}  The following conditions are equivalent.
The first is (\FF.6)  above.

\medskip

(1)\ \ \ $b(x,y_0)$,  $\b(x,y_0)$ and $c(x,y_0)$  vanish on $X$.

\medskip

(2)\ \ \ $H$ is the subset $\{\bll_X\geq0\}$ of $X\times \Jtn$

\medskip

(3)\ \ \ $(\bll f)(x,y_0) = \bll_X(f(x,y_0))$  for all smooth functions $f$ on $Z$.

\medskip

(3)$'$\ \ \ There exists an intrinsic operator $\bll_X'$ on $X$ such that 

\ \ \ \ \ \ \ 
 $(\bll f)(x,y_0) = \bll_X'(f(x,y_0))$  for all smooth functions $f$ on $Z$.

\medskip

(4)\ \ \  $L_{i(x)}\ =\ (i^*)^{-1}(H_x) \ \ \ \forall \ x\in X$.

\medskip

The proof is left to the reader.

\bigskip

\centerline{\bf First Order Restriction}
\medskip

Suppose $F$ is {\sl first order}, that is, $F$ is a subset of $Z\times \JoN$.
By convention the {\sl $F$-subharmonic functions on $Z$} are the same thing as the 
 subharmonic functions for the set  $F\times\Symn\ss \JtN$. 
 If for all compact $K\ss Z$ and $R>0$, 
 $$
 \{(x,r,p) \in F : x\in K, |r|\leq R\}\ \ {\rm is\ compact},
 \eqno{(\FF.8)}
 $$
then $F$ is said to be {\sl coercive}.

If $i:X\hookrightarrow Z$ is defined by $i(x)=(x,y_0)$, and $ H_x \equiv i^*F$
where $i^*$ is restriction of 1-jets, then
$$
 H_x\ =\ \{(r,p) : \exists\, q\ {\rm with\ } (r,(p,q))\in F_{i(x)}\}
\eqno{(\FF.9)}
 $$
If $F$ is coercive, then $H$ is coercive.

\Theorem{\FF.12}  {\sl If $F\ss J^1(Z)$ is coercive and $i:X\hookrightarrow Z$ is defined 
by  $i(x)=(x,y_0)$, then}
$$
u\in F(Z)\qquad\Rightarrow\qquad u\bigr|_X \in  H(X).
$$

\pf The Restriction Hypothesis is easy to verify in this case.  Given $z_0\in X$ and $(r_0,p_0,A_0)$, if 
$$
z_\e\to z_0,\ \ r_\e\to r_0, \ \ {\rm and\ } \ (r_\e, (p_0+A_0(x_\e-x_0), \smfrac 1 \e (y_\e-y_0))\ \in \ F_{z_\e},
$$
then by the coerciveness of $F$ we can extract a subsequence $(z_\e, r_\e, (p_\e,q_\e))$
which converges to  $(z', r', (p',q')) \in F_{z_0}$.  (Here $p_\e \equiv p_0+A_0(x_\e-x_0)$ and 
$q_\e \equiv {1\over \e}(y_\e-y_0)$.)  But $z'=z_0$, $r'=r_0$, and $p'=p_0$.  Hence $(r_0,p_0)\in   H_{x_0}$, which proves the Restriction Hypothesis.\qed

%%%%%%%%%%%%%%%%%%%%%%%%%%%%%%%%%%%%%%%%%%

\def\LQC{\L_q^{\bbc}}

%\vfill\eject
\bigskip

\centerline{\bf  Branches of the Complex Monge-Amp\`ere Equation.}
\medskip

The $q^{\rm th}$ branch $\LQC$of the  complex Monge-Amp\`ere equation
is defined exactly as in the real case  (Examples \BB.5 and  \FF.2) except that the second
derivative $D^2 u$ is replaced by its complex  hermitian part 
${\partial^2 u  \over \partial z_i \partial \bar{z_j}   }  \in {\rm Herm}(\bbc^n)$.
That is, $\LQC$ is the subequation defined by 
$$
\l_q(A_\bbc)\ \geq\ 0 \ \ \ {\rm where}\ \ A_\bbc \ =\ \half(A-JAJ) \ \ {\rm for}\ \ A\in \Sym(\bbr^{2n}).
\eqno{(\FF.10)}
$$
The analogue of (\FF.2) is valid.
$$
\l_q(A_\bbc)\  = \ \inf_W \l_{\rm max} \left( A_\bbc\bigr|_W  \right)
\qquad {\rm for\ all\ }\ A_\bbc \in {\rm Herm}(\bbc^n).
\eqno{(\FF.11)}
$$
where the inf is taken over all complex $q$-dimensional subspaces of $\bbc^n$.
It follows that
$$
{\rm The\ subequation}\ \LQC \ {\rm on}\ \bbc^n \ {\rm restricts\ to
\ the \ subequation}\ \LQC \ {\rm on}\ V 
\eqno{(\FF.12)}
$$
for any complex affine subspace $V$ of dimension $\geq q$.

The characterization Theorem \FF.3 of the $\L_{\rm max}$ subharmonics as the 
``sub'' affine functions, has a natural analogue.  The affine functions are the
solutions to $D^2u=0$.  The {\sl  pluriharmonics}
 are defined to be the solutions of ${  \partial^2 h  \over   \partial z \partial\bar z } =0$.
 Recall that for a simply connected open set $X\ss\bbc^n$
$$
h\ \ {\rm is\ pluriharmonic\ on\ \ } X 
\qquad\iff\qquad 
h\ =\ {\rm Re}F\ \ {\rm with\ } F\ {\rm holomorphic \ on\ \ } X.
\eqno{(\FF.13)}
$$
even when $h$ is only assumed to be a distribution solution.

\Def{\FF.13} A function $u\in\USC(X)$ is {\bf sub-the-pluriharmonics on } $X$
if for each compact subset $K\ss X$ and each  pluriharmonic function $h$ on $X$,
$$
u\ \leq \ h\quad{\rm on}\ \ \partial K
\qquad\iff\qquad 
u\ \leq \ h\quad{\rm on}\ \ K
\eqno{(\FF.14)}
$$

\Prop{\FF.14}  {\sl
A function  $u\in\USC(X)$ with $X\ss\bbc^n$ is $\L^\bbc_{\rm max}$-subharmonic
\ \ $\iff$\ \ $u$ is sub-the-pluriharmonics.
}

\pf
If $u$ is not $\L^\bbc_{\rm max}$-subharmonic on $X$, then it follows from Lemma \BB.2
that  there exist $z_0\in X$, a holomorphic polynomial $F$ of degree 2, with $u(z_0)= {\rm Re} F(z_0)$,
$A\in $ Herm$(\bbc^n)$ with $A<0$ such that
$$
u(z) \ <\ {\rm Re} F(z) + (A(z-z_0), z-z_0) \ \ \ {\rm for}\ z\ {\rm near}\ z_0.
\eqno{(\FF.15)}
$$
Thus $u$ is not sub-the-pluriharmonic Re$F$ on a small ball about $z_0$.
This proves that if $u\in\USC(X)$ is locally sub-the-quadratic-pluriharmonics, 
then $u$ is $\L^\bbc_{\rm max}$-subharmonic.

Now suppose that $u$ is not  sub-the-pluriharmonics on $X$.  That is, for some compact
$K\ss X$ and pluriharmonic function $h$ on $X$, we have
$$
u\ \leq\ h \ \ {\rm on} \ \partial K\quad{\rm but}\quad \sup_K(u-h)\  > \ 0.
$$
This remains true with $h$ replaced by $h-\e|z|^2$ if $\e$ is small enough.
Suppose $z_0$ is a maximum point for  $u-(h-\e|z|^2)$ on $K$.
Adjusting $u$  by subtracting the maximum value at $z_0$, we have 
$u-(h-\e|z|^2) \leq 0$ on $K$ and equal to 0 at $z_0$.  Hence,
$\vf \equiv h-\e|z|^2$ is a test function for $u$ at  $z_0$.
However, since ${\partial^2 \vf  \over \partial z \partial \bar{z}   }(z_0) = -2\e I$,
$u$ is not   $\L^\bbc_{\rm max}$-subharmonic on $X$.  
This proves that if $u$ is  $\L^\bbc_{\rm max}$-subharmonic on $X$, then 
$u$ is sub-the-pluriharmonics on $X$.\qed

\Remark{\FF.15}  The proof shows that the following are equivalent:
\smallskip
(1)\ \ $u$ \  is locally sub-the-quadratic-pluriharmonics,

\smallskip
(2)\ \ $u$ \  is $\L^\bbc_{\rm max}$-subharmonic,

\smallskip
(3)\ \ $u$ \  is  sub-the-pluriharmonics,
\smallskip

\noindent
since (3) \ $\Rightarrow$\ (1) is trivial and we have shown that 
 (1) \ $\Rightarrow$\ (2) \ $\Rightarrow$\ (3).
 \medskip
 
 Combining the Restriction Theorem \FF.1 with the calculation (\FF.12) of the 
 restricted subequation, and with Proposition \FF.14, we have the difficult
 half of the next result.
 
 \Theorem{\FF.16} {\sl
 A function  $u\in\USC(X)$ is $\L^\bbc_{q}$-subharmonic if and only if its
 restriction to each affine complex $q$-plane $V$ is sub-the-pluriharmonics on $X\cap V$.
 }
\pf
Suppose $u$ is not $\L^\bbc_{q}$-subharmonic on $X$.
Then applying Lemma \BB.2 we have (\FF.15) is true with $\l_q^\bbc(A) <0$.
Hence, taking $V$ equal to the span of the first $q$ eigenvectors we see that
$A\bigr|_V<0$, and so $u\bigr|_V$ is not sub-the-pluriharmonics on $V$.\qed

\medskip

Similar results hold for branches of the quaternionic Monge-Amp\`ere Equation.
The details are omitted.

\vskip.3in

%%%%%%%%%%%%%%%%%%%%%%%%%%%%%%%%%%%%%%%%%%%%%%%%
%%%%%%%%%%%%%%%%%%%%%%%%%%%%%%%%%%%%%%%%%%%%%%%%
%%%%%%%%%%%%%%%%%%%%%%%%%%%%%%%%%%%%%%%%%%%%%%%%
%%%%%%%%%%%%%%%%%%%%%%%%%%%%%%%%%%%%%%%%%%%%%%%%
%%%%%%%%%%%%%%%%%%%%%%%%%%%%%%%%%%%%%%%%%%%%%%%%

\vfill\eject

\noindent{\headfont \HH.  The Geometric Restriction Theorem.}

\bigskip

In this section we extend our geometric cases, Theorems \CC.2 and \FF.6, to a full level
of generality.  This is done in three stages delineated as subsections. In the first we treat
restriction to minimal submanifolds of $\rn$ whose tangent planes lie in $\GG$. 
In the next subsection this result is extended to  riemannian manifolds.
 In the final, most general, case, $X\ss Z$ is a $k$-dimensional submanifold 
 of a riemannian manifold $Z$ and the subequation $F = F_\GG\ss J^2(Z)$
 is determined by  a closed subset $\GG\ss G(p,TX)$ of the
bundle of tangent $p$-planes on $Z$ where $p\geq k$.

In all of these cases, because of the additional hypotheses imposed on $X$,
the restricted subequation is also geometrically determined, in fact,  by the set 
$\GG(TX)$ of $\GG$-planes tangent to $X$.  This follows from the algebraic result
Theorem B.3 in Appendix B.

\vskip.3in

\centerline
{\bf Restriction to Minimal $\GG$-Submanifolds.}

\bigskip

  In this subsection the Restriction Theorem \CC.2
will be generalized in two ways.

First , the ``coefficients'' of the subequation are allowed to ``vary''.  That is, a closed subset
$\GG\ss X\times G(p,\rn)$ is given with fibres $\GG_x\ss G(p,\rn)$ defined on an open
set $X\ss\rn$.  Then  the subequation $F$
with fibres $F_x$ is defined by the condition
$$
\trace\left ( A\bigr|_W \right)\ \geq\ 0
\fa W\in \GG_x.
\eqno{(\CC.1)}'
$$
As before, we say that $F$ is {\bf geometrically determined} by 
$\GG\ss X\times G(p,\rn)$.

Second, the affine $\GG$-planes in the RestrictionTheorem  \CC.2 are replaced by 
$\GG$-submanifolds with mean curvature zero.

\Def{\HH.1} A $p$-dimensional submanifold $M$ of $X\ss\rn$ is a {\bf $\GG$-submanifold} if 
$T_xM\in \GG_x$ for each $x\in M$.

\Theorem{\HH.2}  {\sl Suppose $u$ is a $\GG$-plurisubharmonic function on $X\ss\rn$ and $M$ is a 
$\GG$-submanifold of $X$ which is minimal.  Further assume that $\GG\ss X\times G(p,\rn)$ 
has a smooth neighborhood retract which preserves the fibres $\{x\}\times G(p,\rn)$.
If $u$ is $\GG$-plurisubharmonic on $X$, then $u\bigr|_M$ is $\D_M$-subharmonic,
where $\D_M$ is the Laplace-Beltrami operator for the induced metric on $M$.}

\pf
The conclusion is local.  Choose a local orthonormal frame field $e_1,...,e_p$ 
on $M$ and extend it to an orthonormal frame field $e_1,...,e_p$ in a neighborhood $U$
in $\rn$.  Define
$$
W(x)\ =\ \rho(\span\{e_1(x),...,e_p(x)\})
$$
where $\rho$ is the neighborhood retract onto $\GG$.  Then $W(x)$ defines a linear operator
$$
(\bll f )(x) \ \equiv\ \langle P_{W(x)}, \Hess_x f\rangle, \quad{\rm for\ \ } f\in C^\infty(U)
\eqno{(\HH.1)}
$$
(where $P_W$ denotes orthogonal projection onto $W$).  Since  each $W(x)\in \GG$,
we see that if $f$ is $\GG$-plurisubharmonic, then $f$ is $\bll$-subharmonic.
Since $W(x)=T_xM$ for all $x\in M$ we have 
$$
(\bll f)(x) \ =\  \langle T_xM, \Hess_x f\rangle\ =\ (\D_Mf)(x) + (H_M f)(x) \qquad\forall\ x\in M
$$
where $H_M$ is the mean curvature vector field of $M$ (see [HL$_2$] for example). 
Since $M$ is a minimal submanifold, this proves that
$$
(\bll f)(x) \ =\ (\D_Mf)(x)\qquad\forall\ x\in M \ {\rm and
\ } f\in C^\infty(U).
\eqno{(\HH.2)}
$$
Now make a  coordinate change so that $M$ becomes
 $X=\bbr^p\times\{0\}\ss\bbr^p\times\bbr^{n-p}$. By (3)$'$ in Remark \FF.11
 the linear restriction hypothesis is satisfied.  Therefore Theorem \FF.10 implies that
 if an u.s.c. function $u$ is $\GG$-psh, then  $u\bigr|_M$ is $\D_M$-subharmonic.
 \qed

\Remark{\HH.3}  Here we used the obvious fact that  $F_1\ss F_2 \ \ \Rightarrow \ \ 
F_1(X)\ss F_2(X)$  to conclude that if $u$ is $\GG$-plurisubharmonic, then $u$ is $\bll$-subharmonic.

%%%%%%%%%%%%%%%%%%%%%%%%%%%%%%%%%%%%%%%%%%%%%
%%%%%%%%%%%%%%%%%%%%%%%%%%%%%%%%%%%%%%%%%%%%%

\vskip.3in

\centerline
{\bf Riemannian Manioflds.}

\bigskip

The result of the last subsection can be carried over to a completely general version of Theorem \CC.2.
Let $Z$ be a riemannian manifold of dimension $n$ and  $\GG\ss G(p, TZ)$ 
a closed subset of the bundle of tangent $p$-planes on $Z$. 
We again assume   that $\GG\ss G(p, TZ)$ 
admits a smooth neighborhood retraction which preserves the fibres 
of the projection $G(p, TZ) \to Z$.
As before $\GG$ determines a 
natural subequation $F_\GG$ on $Z$ defined by the condition that
$$
\trace\left\{\Hess \, u\bigr|_W \right\} \ \geq \ 0 \fa W \in \GG.
$$
where $\Hess\, u$ denotes the riemannian hessian of $u$.
(See [HL$_{2,6}$] for examples and details.) The corresponding  $F$-subharmonic functions
are again called $\GG${\bf -plurisubhamronic functions}.

By a {\bf $\GG$-submanifold} of $Z$ we mean a $p$-dimensional submanifold
$X\ss Z$ such that $T_xX \in \GG$ for all $x\in X$.
The following result generalizes a basic theorem   in  [HL$_5$]$^*$
\footnote{}{ $ {} \sp{ *}{\rm where} \ F$  was denoted by $\cp^+(\GG)$. } 
for $C^2$-functions to
general upper semi-continuous $\GG$-plurisubharmonic functions.

\Theorem{\HH.4} {\sl Let $X\ss Z$ be a $\GG$-submanifold which is minimal (mean curvature zero).
Then restriction to $X$ holds for  $F_\GG$.  In other words, the restriction of any
 $\GG$-plurisubharmonic function to $X$ is subharmonic in the induced riemannian  metric
 on $X$.}

\pf
Choose local coordinates $z=(x,y)$ on a neighborhood of a fixed
point $(x_0,y_0)$ in $\bbr^p\times \bbr^q$, with $q=n-p$, so that $X$ corresponds
locally to the affine subspace $\{y=y_0\}$.   Choose    a local extension 
of the $\GG$-plane field $TX$ to a $\GG$-plane field $P$ 
defined on a neighborhood $U$ of $(x_0,y_0)$ by taking any  local extension
and composing it with  the neighborhood retraction to $\GG$ as in the proof of Theorem \HH.2.
Consider the linear operator
$$
\bll (u)\ \equiv \ \trace\left\{\Hess \, u\bigr|_P \right\} 
$$
and note that any function which is $\GG$-psh is also $\bll$-subharmonic on $U$.
It will suffice to establish the linear restriction hypothesis for $\bll$.

To see this we note that at points of $X$ the operator $\bll$ can be written as
$$
\bll(u)\ =\ \sum_{i,j=1}^p g^{ij}\left\{   {\partial^2 u\over \partial x_i\partial x_j }
-\sum_{k=1}^p \Gamma_{ij}^k   {\partial u\over \partial x_k}  \right\}
-  \sum_{\a=1}^q  \sum_{i,j=1}^p g^{ij} 
 \Gamma_{ij}^\a   {\partial u\over \partial y_\a}   
\eqno{(\HH.3)}
$$
where $g^{ij}$ denotes the inverse metric tensor and $\Gamma_{ij}^k$ the Christoffel symbols
of the riemannian metric in these coordinates. Equation (\HH.3) can be rewritten as
$$
\bll(u) \ =\ \Delta_X u - H\cdot u
$$where $\Delta_X$ is the Laplace-Beltrami operator for the induced metric on $X$
and $H$ is the mean curvature vector field of $X$.  Since $H\equiv 0$ by hypothesis,
the linear restriction hypothesis (Remark \FF.11 (3)$'$) is satisfied and Theorem \FF.10
 applies to complete the proof.\qed

%%%%%%%%%%%%%%%%%%%%%%%%%%%%%%%%%%%%%%%%%%%%%
%%%%%%%%%%%%%%%%%%%%%%%%%%%%%%%%%%%%%%%%%%%%%

\vskip.3in

\centerline
{\bf The General Geometric Restriction Theorem.}

\bigskip

Just as Theorem \CC.2 was generalized in the last two subsections, the more general result
Theorem \FF.2 can be expanded.
The results of the last two sections can be expanded to a more general situation.
 Let $Z$ and $\GG\ss G(p,TZ)$   be as in the previous subsection.
Fix a submanifold $X\ss Z$
of dimension $m\geq p$ and consider 
the compact subset $\GG(TX) = \{W \in\GG : W  \ss TX\} \ss G(p,TX)$ of
{\sl  $\GG$-planes tangent to $X$}. 
 We say that $X$ is {\bf $\GG$-regular} if each tangent $\GG$-plane at a point $x$
can be extended to a tangent $\GG$-plane field  in a neighborhood of  $x$ in $X$.

The set $\GG(TX)$ defines a subequation $F_{\GG(TX)}$ on $X$ by the requirement that
$$
\trace\left\{\Hess_X \, u \bigr|_W \right\} \ \geq 0 \fa W \in \GG(TX)
$$
for $C^2$-functions $u$, where as before, $\Hess_X$ denotes the riemannian hessian on $X$.

Recall that the {\sl second fundamental form}  $B$ of $X$ is a symmetric bilinear form
on $TX$ with values in the normal bundle $NX$ defined by $B_{V,W} = (\nabla_V \wt W)^N$
where $\wt W$ is any extension of $W$ to a vector field tangent to $X$ (cf. [L]) 
For  $V,W\in T_xX$ the ambient $Z$-hessian and the intrinsic $X$-hessian 
differ by the second fundamental form (cf. [HL$_{2,6}$]), i.e., 
$$
(\Hess_Z  u) (V,W) \ =\  (\Hess_X u) (V,W) + B_{V,W} \, u 
\eqno{(\HH.4)}
$$

\Def{\HH.5} The submanifold $X$ is said to be {\bf $\GG$-flat} if it is $\GG$-regular and
$$
\trace\left\{B\bigr|_W \right\} \ = 0 \fa W \in \GG(TX).
$$

\Theorem{\HH.6. (The Geometric Restriction Theorem)} {\sl
Let $X\ss Z$ be a $\GG$-flat submanifold.  Then the restriction
of any $\GG$-plurisubharmonic function to $X$ is $\GG(TX)$-plurisubharmonic.}

\medskip
\noindent
{\bf Note.}  The simplest interesting case  occurs
when $\dim (X)=p$ and $X$ is a $\GG$-manifold.
Then $X$ is $\GG$-flat if and only if it is minimal   ($\GG$-regularity 
holds automatically).
Thus Theorem \HH.6 generalizes Theorem \HH.4, which in turn contains Theorem \HH.2

Perhaps the next interesting case is that of a real hypersurface $X$ in a K\"ahler manifold
$Z$, where the subset $\GG \ss G_\bbr(2, TZ)$ consists of the complex tangent lines.  We leave it 
to the reader to verify that in this case: {\sl
$X$ is $\GG$-flat if and only if $X$ is Levi-flat}.
\medskip

\noindent
{\bf Proof of Theorem \HH.6.}  From the $\GG$ regularity of $X$ we have the following elementary fact.

\Lemma{\HH.7} {\sl
A function $u\in \USC(X)$ is  $\GG(TX)$psh if and only if for each tangent 
$\GG$-plane field  $W$ defined on an open subset $U\ss X$, the function $u\bigr|_U$ is
$\bll_W$-subharmonic, where $\bll_W$ is the linear subequation on $U$defined by
$\bll_W (v) \equiv \tr_W \{\Hess_X  v \}            \geq 0$ for $v\in C^2$.
}
\pf
($\Leftarrow$)
Let $\vf$ be a test function for $u$ at $x_0\in X$. Fix $W_0\in \GG(T_{x_0}X)$.  
Extend $W_0$ to a local $\GG(TX)$-plane field $W$.
Then by assumption  $\tr_{W }\{\Hess_X \vf\} \geq 0$. This proves that $\tr_{W_0}\{\Hess_X \vf\} \geq 0$
for all $W_0\in \GG(T_{x_0}X)$, i.e., $\Hess_{x_0} \vf \in F_{\GG(T_{x_0}X)}$.   

($\Rightarrow$)
Suppose $u$ is $\GG(TX)$-psh, and let $W$ be a tangent $\GG$ plane 
field defined on an open set $U\ss X$. Fix $x_0\in U$ and choose a test function $\vf$ for
$u$ at $x_0$.  Since $u$ is $\GG(TX)$-psh, we have $\tr_{W_0}\{\Hess_X \vf\} \geq 0$
for all $W_0\in \GG(T_{x_0}X)$.
Hence $u$ is $\bll_W$-subharmonic on $U$.   \qed
\medskip

The remainder of the proof of Theorem \HH.6 now closely follows the argument given
for the proof of Theorem \HH.4, by choosing similar coordinates and extending the intrinsic
operators $\bll_W$ into $Z$.\qed

\Ex{\HH.8. ($\GG$-regularity is necessary)} Let $\GG = \{x$-axis$\}$ in $\bbr^2$, and set 
$X =\{(x,y) : y=x^4\}$.  Then $X$ has a tangent $\GG$-plane only at the origin. 
The second fundamental form (i.e., the curvature) is zero at the origin, however
 $\GG$-regularity clearly fails.  Restriction also fails. Consider the strictly $\GG$-psh
function $u(x,y) = \e x^2 -{|y|}^{1\over2}$.  Then $u\bigr|_X = u(x,x^4) = -(1-\e)x^2$ in the parameter
$x$, and one sees easily that for $\e$ small, $\Hess_0\, u = {d^2u\over ds^2}(0)<0$
(where $s = $ arc-length parameter).

%\vfill\eject
\vskip .3in

%%%%%%%%%%%%%%%%%%%%%%%%%%%%%%%%%%%%%%%%%%%%%
%%%%%%%%%%%%%%%%%%%%%%%%%%%%%%%%%%%%%%%%%%%%%
%%%%%%%%%%%%%%%%%%%%%%%%%%%%%%%%%%%%%%%%%%%%%
%%%%%%%%%%%%%%%%%%%%%%%%%%%%%%%%%%%%%%%%%%%%%

\noindent{\headfont \LL.\   Jet Equivalence of Subequations.}

\medskip
In this section and the next one we suppose that a subequation $F$ is given on a smooth manifold $Z$.
No riemannian assumption will be made.
In particular, $F$ is a closed subset of the 2-jet bundle $J^2(Z)$.  The 0-jet bundle $\bbr$ splits off
as $J^2(Z) = \bbr \oplus J^2_{\rm red}(Z)$ leaving the bundle of reduced 2-jets 
$J^2_{\rm red}(Z)$.  The bundle of reduced 1-jets is simply $T^*Z$ the cotangent bundle of $Z$.

\bigskip
\centerline{\bf    Restriction}
\medskip

If $X$ is a submanifold of $Z$, let $i^*_X$ denote the restriction of 2-jets   to $X\ss Z$. 
Then
$$
\matrix
{
0  &  \arr &  \Sym(T^*Z) &     \arr &  J^2_{\rm red}(Z)    &     \arr &      T^*Z
&     \arr &    0   \cr         
\ \cr   
\  &   \ &  \downarrow i^*_X &     \ &   \downarrow i^*_X    &     &      \downarrow i^*_X
&     &       \cr            
\ \cr
0  &  \arr &  \Sym(T^*X) &     \arr &  J^2_{\rm red}(X)    &     \arr &      T^*X
&     \arr &    0   \cr            
}
\eqno{(\LL.1)} 
$$
is commutative with exact rows.  Note that $i^*_X : \Sym(T^*Z) \to \Sym(T^*Z)$ is  the restriction of quadratic forms on $TZ$ to quadratic forms on $TX$, and that the quotient map 
$i^*_X : T^*Z \to T^*X$ is restriction of 1-forms.

\bigskip
\centerline{\bf    Automorphisms}
\medskip
 To begin, an automorphism of the jet bundle $J^2(Z) = \bbr\oplus J^2_{\rm red}(Z)$
  is required to split as the identity on the 0-jet factor $\bbr$ and an automorphism of the reduced
  jet bundle $J^2_{\rm red}(Z)$.  Hence it suffices to define automorphisms of the reduced jet bundle.
  
  \Def{\LL.1}  An {\bf automorphism}  of $J^2_{\rm red}(Z)$ is a bundle isomorphism
  $\Phi: J^2_{\rm red}(Z) \to J^2_{\rm red}(Z)$ which maps the 
  subbundle $\Sym(T^*Z)$ to itself and has the further property that this restricted
  isomorphism $\Phi: \Sym(T^*Z)\to\Sym(T^*Z)$ is induced by a bundle isomorphism
  $$
  h=h_{\Phi} : T^*Z\ \arr\ T^*Z.
  \eqno{(\LL.2)}
  $$
  \smallskip
  
This means that for $A\in \Sym(T^*Z)$,
 $$
  \Phi(A) \ =\ hAh^t,
  \eqno{(\LL.3)}
  $$
that is,
$$
\Phi(A)(v,w)  \ =\ A(h^tv, h^tw) \qquad {\rm for }\ \ v,w \in TZ.
$$
 
 \medskip
 
 Because of the upper short exact sequence in (\LL.1)
 each automorphism $\Phi$ of $J^2_{\rm red}(Z)$ induces a bundle isomorphism
  $$
 g=g_\Phi: T^*Z\ \to\ T^*Z.
  \eqno{(\LL.4)}
  $$
 This bundle isomorphism is {\sl not} required to agree with $h$ in (\LL.2).

\Lemma{\LL.2} {\sl  The automorphisms of $J^2(Z)$ form a group. They are the sections 
of the bundle of groups whose fibre at $z\in Z$ is the group of automorphisms
of $J_z^2(Z)$ defined   above.}

\pf
 See [HL$_6$, \S 4].

\Prop{\LL.3} {\sl With respect  to any  splitting 
$$
J^2(Z) \ =\ \bbr\oplus T^*Z \oplus \Sym(T^*Z)
$$
of the upper short exact sequence (\LL.1), a bundle automorphism has the form
$$
\Phi(r, p,  A) \ =\ (r, gp, hAh^t + L(p))
\eqno{(\LL.5)}
$$
where $g$  and $h$ are smooth sections of the bundle $\End(T^*Z)$ and $L$ is a smooth section of the bundle 
$\Hom (T^*Z, \Sym(T^*Z))$. }
\pf Obvious.

\Ex{1} The trivial 2-jet bundle on $\rn$ has fibre
$$
 \bbj^2  = \bbr\times\rn\times \Symn.
 $$
with automorphism group
 $$
{\rm Aut}(\bbj^2) \ \equiv\  {\rm GL}_n\times {\rm GL}_n\times \Hom(\rn, \Symn)
 $$
 where the  action is  given by
 $$
 \Phi_{(g,h,L)}(r,p,A) \ =\ (r, \ gp, \ hAh^t +L(p)).
 $$
 and the group law is
 $$
(\bar g, \bar h, \bar L)\cdot  (g,h,L)\ =\ (\bar g g,\  \bar h h, \bar h L \bar h^t +\overline{ L}\circ g)
 $$

\Ex{2}  Given a local coordinate system $(x_1,...,x_n)$ on an open set
$U\ss Z$, the  canonical trivialization 
$$
J^2(U)\ =\ U\times\bbr\times\rn\times\Symn
\eqno{(\LL.6)}
$$
is determined by the coordinate 2-jet $J_xu  = (u,Du,D^2u)$ evaluated at $x$. 
With respect to this splitting, every automorphism is of the form
$$
\Phi(u, Du, D^2u) \ =\ (u, \ g Du, \ h\cdot D^2u\cdot h^t + L(Du))
\eqno{(\LL.7)}
$$
where $g_x, h_x \in {\rm GL}_n$ and $L_x : \rn\to \Symn$ is linear for each point $x\in U$.

\bigskip
\centerline{\bf Jet Equivalence}
%\medskip

\Def{\LL.4}  Two subequations $F, F'\ss J^2(Z)$ are {\bf  jet equivalent} if 
there exists an automorphism  $\Phi:J^2 (Z)\to  J^2(Z)$ with 
$\Phi(F)=F'$.

\Def{\LL.5}  A subequation $F\ss J^2(Z)$ is {\bf locally jet equivalent to a constant coefficient
subequation} if each point $x$ has a distinguished coordinate neighborhood $U$ so that 
$F\bigr|_U$ is jet equivalent to a constant coefficient subequation $U\times \bbf$ in those distinguished
coordinates.

\Lemma{\LL.6} {\sl  Suppose $Z$ is connected and $F\ss J^2(Z)$ is  locally jet equivalent to a constant coefficient
subequation. Then there is a subequation $\bbf\ss \bbj^2$, unique up to equivalence,
such that $F$ is locally jet equivalent to $U\times \bbf$ on every distinguished coordinate chart.}

\pf 
In the overlap of any two distinguished charts $U_1\cap U_2$ choose a point $x$.  Then the 
local equivalences $\Phi_1$ and $\Phi_2$, restricted to $F_x$, determine an equivalence 
from $\bbf_1$ to $\bbf_2$.  Thus the local constant coefficient equations on these
charts are all equivalent, and they can be made equal by applying the appropriate
constant equivalence on each chart.\qed

\Remark{\LL.7}  The notion of jet equivalence arises naturally when considering
the group of germs  of diffeomorphisms  which fix  a point $x_0$,  acting on $J_{x_0}^2$.
 Namely, if $\vf$ is a local diffeomorphism
fixing $x_0$, then in local coordinates (as in Example 2 above) the right action on $J_{x_0}^2$, induced
by the pull-back $\vf^*$ on 2-jets,  is given by (\LL.7) where $g_{x_0}=h_{x_0}$ is the transpose of the
Jacobian matrix $(({\partial \vf^i\over \partial x_j}))$
and $L_{x_0}(Du) = \sum_{k=1}^n  {\partial^2 \vf^k\over \partial x_i\partial x_j}(x_0) u_k$.
Thus  with jet coordinates $(r,p,A)$ at $x_0$
$$
\vf^*(r,p,A) \ =\ (r, \ gp, \ gAg^t +   D^2_{x_0}(\vf) \cdot p).
$$
Note, however, that this applies {\sl only} at the fixed point $x_0$.
%A jet equivalence induced by an automorphism on an open set is necessarily
%the identity.

\medskip\noindent
{\bf Cautionary Note.}  A local equivalence $\Phi: F\to F'$ does not take $F$-subharmonic functions
to $F'$-subharmonic functions.   In fact, for  $u \in C^2$,  $\Phi(J^2 u)$ is almost never the 2-jet of a 
function. It happens if and only if  $\Phi(J^2 u) = J^2u$.

%%%%%%%%%%%%%%%%%%%%%%%%%%%%%%%%%%%%%%%%%%%%%
%%%%%%%%%%%%%%%%%%%%%%%%%%%%%%%%%%%%%%%%%%%%%
%%%%%%%%%%%%%%%%%%%%%%%%%%%%%%%%%%%%%%%%%%%%%
%%%%%%%%%%%%%%%%%%%%%%%%%%%%%%%%%%%%%%%%%%%%%

\bigskip

\centerline{\bf  Relative Automorphisms and Relative Jet Equivalence}
\medskip

Suppose now that $i:X\hookrightarrow Z$ is an embedded submanifold.

\Def{\LL.8}  A {\bf relative automorphism of $J^2(Z)$ with respect to $X$} is an automorphism 
$\Phi:J^2(Z)\to J^2(Z)$  such that on $X$ the diagram 
$$
\matrix
{
J^2(Z) & \harr {\Phi}\ & J^2(Z)  \cr
i^* \downarrow & \  & \downarrow i^*  \cr
J^2(X) &  \harr {\vf}\  & J^2(X)  \cr
}
$$
commutes for some automorphsim $\vf:J^2(X)\to J^2(X)$.

\medskip
Relative automorphisms with respect to $X$ are a subgroup of the automorphisms of $J^2(Z)$.

Fix a splitting $J^2(Z) \ =\ \bbr\oplus T^*Z \oplus \Sym(T^*Z)$, and let $g$, 
  $h$ and $L$ be associated to an automorphism $\Phi$ as in Proposition \LL.4.
  Then one easily checks that: {\sl $\Phi$ 
 is a  relative automorphism of $J^2(Z)$ with respect to $X$ if and only if }
 $$
 g^t(TX)\ \ss\ TX, \qquad h^t (TX)\ \ss\ TX 
 \and
 L_{N^*X, \Sym(T^*X)}\ =\ 0.
 \eqno{(\LL.8)}
 $$
Here $L_{N^*X, \Sym(T^*X)}$ denotes the restriction of $L$ to $N^*X$ followed by the 
restriction of quadratic forms in $\Sym(T^*Z)$ to $\Sym(T^*X)$.

\Def{\LL.9} Two subequations $F,F'\ss J^2(Z)$ are {\bf jet equivalent 
modulo $X$} if  $F' = \Phi(F)$ for some relative automorphism $\Phi$ with respect to $X$.
\medskip

If  $F,F'\ss J^2(Z)$ are  jet equivalent 
modulo $X$, then the induced subequations $H=i^*F$ and $H'=i^*F'$ are jet equivalent on $X$.

By an {\sl adapted coordinate neighborhood} of a point $z_0=(x_0, y_0) \in X$ we mean a local coordinate
system $z=(x,y)$ on a neighborhood $U$ of $z_0$ such that $X\cap U
= \{ (x,y) : y=y_0\}$.

\Def{\LL.10}  The  subequation $F\ss J^2(Z)$ is {\bf locally jet equivalent  modulo $X$ to a constant coefficient
subequation} if each point  in $X$ has an adapted coordinate neighborhood $U$ so that 
$F\bigr|_U$ is  jet equivalent  modulo  $X$
to a constant coefficient subequation $U\times \bbf$ in those adapted
coordinates.

\medskip

Now we examine what this means in more detail. 
Suppose that $z=(x,y) \in \bbr^N= \rn\times\bbr^m$ is the adapted coordinate system
and $\Phi : J^2(U)\to J^2(U)$ is the jet equivalence modulo $X$.
By Proposition \LL.3, $\Phi$ acting on a coordinate 2-jet $(u,Du,D^2u)$
must be of the form 
$$
\Phi(J)\ =\ \Phi(u,Du,D^2u)\ =\ (u, gDu, hD^2uh^t + L(Du)).
\eqno{(\LL.8)}
$$
Moreover, we have
$$
J\in F\quad\iff\quad \Phi(J)\in \bbf.
\eqno{(\LL.9)}
$$
With respect to the splitting $\bbr^n\times \bbr^m$ into $x$ and $y$ coordinates,  each
coordinate 2-jet J
can be written as  
$$
J\ =\  \left(r, (p,q), \left(\matrix {A&C\cr  C^t &B}\right)\right),
\and
i^*(J)\ =\ (r,p,A).
$$
is the restriction of $J$ to $X$.
The sections $g$ and  $h$   can be written in block form as
$$
g \ =\ \left(\matrix { g_{11} & g_{12} \cr g_{21} & g_{22}}\right)
\and
h \ =\ \left(\matrix { h_{11} & h_{12} \cr h_{21} & h_{22}}\right).
 \eqno{(\LL.10)}
$$
Also $L$ can be decomposed into the sum $L=L'+L''$ where
$L' \in\End(\rn, \Sym(\bbr^N))$ and 
$L'' \in\End(\bbr^m, \Sym(\bbr^N))$.  Each of  $L,L', L''$ can be blocked into
$(1,1)$, $(1,2)$, $(2,1)$, $(2,2)$ components in $\Sym(\rn\oplus\bbr^m)$, 
in analogy with $g$ and $h$ above.

Now we can compute the restriction $i^*\Phi(J)$ of $\Phi(J)$.  Namely,
$$
\eqalign
{
i^*\Phi(J)\ =\  (r, \ g_{11}p   &+g_{12}q, \ \ h_{11}A h_{11}^t +h_{12}C^t h_{11}^t  \cr
&+ h_{11}C h_{12}^t +h_{12}B h_{12}^t
+L_{11}'(p) + L_{11}''(q)  )  
}
\eqno{(\LL.11)}
$$
In order for $\Phi$ to be a jet equivalence modulo $X$
this must agree with an automorphism $\vf : J^2(U\cap X) \to J^2(U\cap X)$,
which is the case if and only if on $X$
$$
g_{12}\ =\ 0, \quad h_{12}\ =\ 0, \quad {\rm and}\quad L_{11}''\ =\ 0
\eqno{(\LL.12)}
$$
so that
$$
\vf(r,p,A)  \ =\  \left(r, \ g_{11}p,\ h_{11}Ah_{11}^t + L_{11}'(p)\right)
\eqno{(\LL.13)}
$$

\medskip
\noindent
{\bf Final Note. \LL.11. (Affine Jet Equivalence).} The above discussion extends easily to
the more general case of {\sl affine automorphisms}.  The affine automorphism group 
is an extension of   the automorphism 
group of $J^2(Z)$   via bundle translations by sections of $J^2(Z)$. (See [HL$_6$, \S  6.3] for details.)

%%%%%%%%%%%%%%%%%%%%%%%%%%%%%%%%%%%%%%%%%%%%%
%%%%%%%%%%%%%%%%%%%%%%%%%%%%%%%%%%%%%%%%%%%%%
%%%%%%%%%%%%%%%%%%%%%%%%%%%%%%%%%%%%%%%%%%%%%
%%%%%%%%%%%%%%%%%%%%%%%%%%%%%%%%%%%%%%%%%%%%%

\vfill \eject

\noindent{\headfont \MM.\  The  Restriction Theorem for Subequations Derivable from a Euclidean Model.}
\medskip

The next result does not include the Geometric Restriction Theorem \FF.6 since the
subset $\GG\ss G(p,TX)$ may not even be a subbundle of $G(p,TX)$.  However,
it applies to some interesting non-geometric cases, and to some cases of a geometric 
but non-riemannian type.  The non-geometric  application is given in the 
next Section \NN. The non-riemannian application with a geometric flavor is given in the 
separate paper [HL$_8$] where we prove that restriction holds for the intrinsically defined plurisubharmonic functions on an almost complex manifold.

\Theorem{\MM.1}  {\sl Let $i:X\hookrightarrow Z$ be an embedded submanifold and $F\ss J^2(Z)$ 
a subequation.  Assume that 
  $F$ is  locally jet equivalent  modulo $X$ to a constant coefficient
subequation $\bbf$.   Set $\bbh \equiv \overline{ i^*\bbf}$ .   Then 
$H\equiv \overline{ i_X^*F}$ is locally jet equivalent to the constant coefficient subequation $\bbh$, 
and restriction holds.  That is, }
$$
u \ \ {\sl is}\ F\ {\sl subharmonic\ on\ }\ Z\qquad\Rightarrow\qquad 
u\bigr|_X \ \ {\sl is}\ H\ {\sl subharmonic\ on\ }\  X
$$

\pf
Adopt the notation following Definition \LL.10.
By hypothesis (\LL.12) we have that
$$
g_{12}(x,y) \ \ {\rm and}\ \  h_{12}(x,y) \ {\rm are} \ O(|y-y_0|)
\and
{L''}_{11}(x,y) \ =\ O(|y-y_0|).
\eqno{(\MM.1)}
$$

We now show that $F$ satisfies the Restriction Hypothesis.  Fix $(r_0, p_0, A_0) \in J^2_{x_0}(X)$
and suppose there are  sequences $z_\e = (x_\e,y_\e)$ and $r_\e$ with 
$$
J_\e \ =\ \left( r_\e, \left( p_0 + A_0(x_\e-x_0), {y_\e-y_0\over  \e}\right), \left( 
\matrix{A_0 & 0 \cr 0 & {1\over \e} I}
  \right)\right)\ \in \ F_{z_\e}
  \eqno{(\MM.2)}
$$
and 
$$
  x_\e\ \to\  x_0,\ \ {{|y_\e-y_0|^2}\over\e}\ \to\ 0, \ \  r_\e\ \to\ r_0,
    \eqno{(\MM.3)}
$$
as   $\e \to 0$. Now (\MM.2) is equivalent to the fact that
$$
\Phi_{z_\e}(J_\e) \ \in \ \bbf  \fa \e.
$$
This means that the  $(1,1)$-component 
$$
i^* \Phi_{z_\e}(J_\e) \ \in \ i^* \bbf    \fa \e.
  \eqno{(\MM.4)}
$$
To show that $(r_0,p_0,A_0) \in H_{z_0} =  \overline{i^*_X F_{z_0}}$ it will suffice  to show that  
$$
i^* \Phi_{z_\e}(J_\e) \ \   {\rm converges\  to \ \ } \vf(r_0,p_0, A_0)\ \ 
 {\rm as\ } \e\to 0.
  \eqno{(\MM.5)}
$$
Write
$$
i^* \Phi_{z_\e}(J_\e) \ =\ \left(r_\e, p_\e, A_\e   \right).
$$
By (\LL.11) 
$$
p_\e \ =\ g_{11}(z_\e)(p_0+A_0(x_\e-x_0))+ g_{12}(z_\e)\smfrac 1 \e(y_\e-y_0).
$$
Now  (\MM.1) and (\MM.3)  imply that $p_\e \to g_{11}(z_0) p_0$. Furthermore, by (\LL.11) 
$$
\eqalign
{
A_\e\ &=\ h_{11}(z_\e) A_0 h_{11}^t (z_\e) + \smfrac 1 \e h_{12}(z_\e)h_{12}^t(z_\e)  \cr
& \qquad \qquad \qquad 
+{L'}_{11}(z_e)\cdot (p_0+A_0(x_\e-x_0)) +  {L''}_{11}(z_e) \cdot ((\smfrac 1 \e(y_\e-y_0)).\cr
}
$$
Again by (\MM.1) and (\MM.3) we have 
$A_\e \to h_{11}(z_0) A_0 h_{11}^t (z_0) + {L'}_{11}(z_0)\cdot p_0$.
Since $\vf_{z_0}(r_0,p_0, A_0) = (r_0, \ g_{11}(z_0)p_0,  \ h_{11}(z_0) A_0 h_{11}^t (z_0)
 + {L'}_{11}(z_0)\cdot p_0)$, this completes the proof.\qed

%\vfill\eject
\vskip .3in

%%%%%%%%%%%%%%%%%%%%%%%%%%%%%%%%%%%%%%%%%%%%%%%
%%%%%%%%%%%%%%%%%%%%%%%%%%%%%%%%%%%%%%%%%%%%%%%
%%%%%%%%%%%%%%%%%%%%%%%%%%%%%%%%%%%%%%%%%%%%%%%
%%%%%%%%%%%%%%%%%%%%%%%%%%%%%%%%%%%%%%%%%%%%%%%
%%%%%%%%%%%%%%%%%%%%%%%%%%%%%%%%%%%%%%%%%%%%%%%

\noindent{\headfont  \NN.   Applications of this Last  Restriction Theorem.}\medskip

\bigskip

The Second  Restriction Theorem has a number of interesting applications.
One is to the universally defined subequations on manifolds with topological
$G$-structure (as in [HL$_6$]).

We begin with the case of universal riemannian subequations.
By a {\sl euclidean model} we mean a closed subset
$$
\bbf \ \ss \ \bbj^2_N = \bbr\times \bbr^N\times \Sym(\bbr^N)
\eqno{(\NN.1)}
$$
 with the properties that:
\medskip

(1) \ \    $\bbf+ (\bbr_- \times \{0\}\times \cp) \ss\bbf$,  \ \ \ where  $\cp\equiv \{A \in  \Sym(\bbr^N): A\geq0\}$,  
\smallskip

(2)\ \ \  $\bbf = \overline {\Int  \bbf }$, \ \ \ and 
\smallskip

(3) \ \  $\bbf$ is invariant under the natural action of O$_N$ on $\bbj^2_N$.
\medskip

 \noindent
 Let $Z$ be a riemannian manifold of dimension $N$ and recall the canonical
splitting 
$$
J^2(Z) \  =  \ \bbr\times T^*Z\times \Sym(T^*Z)
\eqno{(\NN.2)}
$$ 
given by the riemannian hessian
$$
(\Hess \, u)(V,W)\ \equiv\ VWu -(\nabla_V W)u
\eqno{(\NN.3)}
$$
(for vector fields $V$ and $W$; see [HL$_6$].)

\Def{\NN.1} 
The model subequation $\bbf$ in (\NN.1) is {\bf universal} because it canonically determines
a subequation $F\ss J^2(Z)$ on any riemannian $N$-manifold $Z$ 
by the requirement that 
$$
Ju_z = (u(z), (du)_z, \Hess_z u) \  \in \ F_z \quad \iff\quad [u(z), (du)_z, \Hess_z u]  \ \in  \ \bbf
\eqno{(\NN.4)}
$$
where $[u(z), (du)_z, \Hess_z u]$  denotes the coordinate representation
of  $(u(z), (du)_z, \Hess_z u)$ with respect to any orthonormal basis 
of $T_zZ$.  We call $F$ {\bf the subequation on $Z$ canonically determined by 
$\bbf$}.

\Theorem {\NN.2. (Restriction for Universal Riemannian Subequations)} {\sl
Let $Z$ be a riemannian manifold of dimension $N$ and $F\ss J^2(Z)$ a subequation canonically
determined by an O$_N$-invariant universal subequation $\bbf\ss\bbj^2_N$ as above.
Then restriction holds for $F$ to any totally geodesic submanifold $X\ss Z$.
}

\pf     
The theorem is local, so we may restrict to the case where 
$$
\eqalign
{
Z\ &\equiv \ \{ x= (x',x'') \in \bbr^n\times \bbr^m : |x'| < 1, |x''| < 1\},\ \  {\rm and} \cr
X\ &\equiv \ \{ x= (x',0) \in \bbr^n\times \bbr^m : |x'| < 1\}, \cr
}
$$
with $n+m=N$. 
We may furthermore assume that    
$$
\partial_i' \ \perp\ \partial_j''  \quad {\rm along \ }X \fa i,j
\eqno{(\NN.5)}
$$
 in the given metric on $Z$  where 
  $$
  \partial_i' \equiv  {\partial \over \partial x_i'} \and
 \partial_j'' \equiv  {\partial \over \partial x_j''}.
 $$
To see this we choose our coordinates as follows.  First choose a local coordinate map
 $\vf : \{x', |x'| \leq 1\} \to X$. Fix a basis $\nu_1,...,\nu_m$ 
 of the normal space to $X$ at $\vf(0)$
and extend them to normal vector fields $\nu_1,...,\nu_m$ on $X$ by 
parallel translation along the curves corresponding to rays from the  origin 
in the disk   $\{x', |x'| \leq \ 1\}$. Applying the exponential map to 
 $x_1'' \nu_1\vf((x')) + \cdots +  x_m''\nu_m(\vf(x'))$ gives the desired coordinates for $|x''|  < $ some $\e$.
 (Of course, one can then renormalize  to $|x''| <1$.)

We now choose an orthonormal frame field $(e_1,...,e_{n+m}) = (e_1',...,e_n', e_1'',...,e_m'')$
on $Z$ (with respect to the given metric) so that along $X$
$$
e_1',...,e_n' \ \ {\rm are \ tangent\ to\ \ } X\and e_1'',...,e_m'' \ \ {\rm are \ normal\ to\ \ } X.
\eqno{(\NN.6)}
$$
Our subequation  $F\ss J^2(Z)$ is then given explicitly by the condition
$$
\bigl(u, (e_1u, ..., e_{n+m}u),  \Hess\, u(e_i,e_j) \bigr)_z \ \in\ \bbf
\eqno{(\NN.7)}
$$
for $z\in Z$.
We now write
$$
e_i\ =\ \sum_{j=1}^{n+m} h_{ij} \partial_j \qquad{\rm for}\ \  \  i=1,...,n+m
$$
where $\partial \equiv (\partial', \partial'')$. From (\NN.5) we have that the matrix $h$ decomposes
as
$$
h\ =\ \left(   
\matrix
{
h' & 0\cr
0 & h''\cr
}
\right)
\qquad {\rm along\ \ } X.
\eqno{(\NN.8)}
$$
We now compute that
$$
e_i u \ =\ \sum_j h_{ij} \partial_j u, \qquad{\rm and}
$$
$$
\eqalign
{
(\Hess\, u)(e_i, e_j) \ &= \ (\Hess\, u)\left(\sum_k h_{ik} \partial_k , \sum_\ell h_{j\ell} \partial_\ell \right) 
\ =\ \sum_{k,\ell} h_{ik} h_{j\ell} (\Hess\, u)(\partial_k, \partial_\ell)   \cr
&= \  \sum_{k,\ell} h_{ik} h_{j\ell} \left\{  \partial_k\partial_\ell u - \left(\nabla_{\partial_k} \partial_\ell\right) u\right\}\cr
&= \  \sum_{k,\ell} h_{ik} h_{j\ell} \left\{  \partial_k\partial_\ell u - \sum_m \G_{k\ell}^m \partial_m u\right\}\cr
}.
$$
where $\G = \{\G_{k\ell}^m\}$ are the classical Christoffel symbols.  Expressed briefly, we have that
$$
e\cdot u  \ =\ h Du
\and
(\Hess\, u)(e_*, e_*) \ =\ h(D^2 u) h^t - \wt{\G} \cdot Du
$$
where $\wt \G \equiv h \G h^t$. Thus our condition (\NN.7) can be rewritten in terms of the
coordinate jets as 
$$
\left(u, \, h Du, \, h(D^2u)h^t - \wt \G\cdot Du \right) \ \in\ \bbf.
\eqno{(\NN.9)}
$$
This says precisely that our subequation $F$ is jet equivalent to the  constant coefficient
subequation $\bbf$ in these coordinates.  

We claim that this is an equivalence mod $X$.  For this we must establish the conditions in (\LL.12).
Note first that in this case $g=h$ and $h_{12} =0$ by (\NN.8). For the last condition we use the fact
that $X$ is totally geodesic.  This means precisely  that 
$$
\nabla_{\partial_i'}  \partial_j' \ =\ \sum_{k=1}^n \G_{ij}^k \partial_k' \qquad{\rm along\ \ } X,
$$
i.e. $\nabla_{\partial_i'}  \partial_j' $ has no normal components along $X$ for all $1\leq i,j\leq n$.
  This is exactly the third condition in (\LL.12).
  
Theorem \NN.2 now follows from Theorem \MM.1.
\qed

\medskip

Theorem \NN.2 can be extended to the case where the riemannian manifold $Z$ has a
topological reduction of the structure group to a subgroup 
$$
G\ \ss\ {\rm O}_N.
$$
Such a reduction consists of an open covering $\{U_\a\}_\a$ of $Z$ and an orthonormal
tangent frame field $e^\a = (e^\a_1,...,e_N^\a)$ given on each open set $U_\a$ with the property
that the change of framings
$$
g_{\a\b} : U_\a \cap U_\b \ \arr\ G\ \ss\ {\rm O}_N
$$
take their values in $G$. 

The local frame fields $e_a$ are called {\bf admissible}.  Note that if $e$ on $U$ is an admissible
frame field, one can add to the family of admissible framings, any frame field of the form
$ge$ where $g:U\to G$ is a smooth   map. We assume that our $G$-structure has
a maximal family of admissible frame fields.

\Def{\NN.3} Suppose $Z$ has a topological $G$-structure.  A submanifold $X\ss Z$ 
is called {\bf $G$-adaptable}  if for every point $z\in X$ there is an admissible framing
$e$ on a neighborhood $U$ of $z$ such that on $X\cap U$
$$
e_1,...,e_n \ \ {\rm are \ tangent\ to\ \ } X\cap U \and e_{n+1},...,e_N \ \ {\rm are \ normal\ to\ \ } X\cap U.
\eqno{(\NN.10)}
$$

\Ex{\NN.4}  Suppose $G = {\rm U}_m \ss {\rm O}_{2m}$.  Having a ${\rm U}_m$-structure on $Z$
is equivalent to having an orthogonal almost complex structure $J:TZ\to TZ$, $J^2\equiv -I$ on $Z$.
A ${\rm U}_m$-adaptable submanifold $X\ss Z$  is simply an almost complex submanifold,
i.e., having the property that $J(T_xX)=T_xX$ for all $x\in X$.

\medskip

On a manifold with topological $G$-structure, we can enlarge the set of 
universal subequations by replacing property (3) above with 
\smallskip
(3)$'$ \ \  $\bbf$ is invariant under the natural restricted action of  $G$ on $\bbj^2_N$.
\smallskip
As above any such set $\bbf$ determines a subequation $F$ on $Z$.

\Theorem{\NN.5}  {\sl
Let $Z$ be a riemannian manifold with topological $G$-structure, 
and $F\ss J^2(Z)$ a subequation canonically
determined by a  $G$-invariant universal subequation $\bbf\ss\bbj^2_N$
satisfying (1), (2) and (3)$'$. 
Then restriction holds for $F$ to any totally geodesic $G$-adaptable submanifold $X\ss Z$.
}
 \pf The proof exactly follows the one given for Theorem \NN.2. One merely has to choose
 the local frame field $e$ with property (\NN.6) to be an admissible field (cf. (\NN.10)).
 Details are left to the interested reader.\qed

\Note{\NN.6}  Every almost complex manifold $(Z,J)$ admits many almost complex submanifolds
of dimension one (pseudo-holomorphic curves) by a classical result of Nijenhuis and Woolf [NW].
In fact there exist  pseudo-holomorphic curves in every complex tangent direction at every point.
It is standard to define an upper semi-continuous function to be plurisubharmonic if its
restriction to every such curve is subharmonic. Using Theorem 8.1 above, the authors have
proved  in  [HL$_8$] that this standard definition of plurisubharmonicity
coincides with the viscosity definition coming from 
an intrinsically defined subequation $F(J)$ on $Z$.  They also show in  [HL$_8$] that the standard
plurisubharmonic functions are, in a precise sense, equivalent to the plurisubharmonic
distributions on $(Z,J)$.

 \medskip
 
 Theorem \NN.2 asserts that every universal riemannian subequation
 satisfies restriction to totally geodesic submanifolds. Of course 
 if the submanifold $X$ is too small,
 this restriction is trivial, i.e., $i^*F = J^2(X)$.  One extreme example of this is the Laplace-Beltrami 
 equation given by $\bbf = \{(r,p, A) : \tr A\geq0\}$ where all submanifolds (even hypersurfaces) are too small. Nevertheless, there are also many subequations
 which have interesting restrictions.  One such is the classical $\bbf = \{(r,p, A) : A\geq0\}$
 corresponding to riemannian convex functions.  This branch of Monge-Amp\`ere falls under the aegis of Geometric Restriction Theorem \JJ.2, but the other branches are not covered by previous results.
 Recall the constant coefficient case Example \BB.5/\FF.2.

 \Ex{\NN.7. (The Monge-Amp\`ere Equation)}  Given $A\in \Sym(\bbr^N)$, let $\l_1(A) \leq \cdots \leq \l_N(A)$ denote as before 
the {\sl  ordered eigenvalues} of $A$.  Define for $\mu\in \bbr$
$$
{\bf \Lambda}^{\mu}_q \ \equiv\ \{(r,p, A) \in \bbj^2_N : \l_q(A) \geq \mu\}.
$$
Let $\Lambda^{\mu}_q (Z)$ be the induced subequation on the riemannian manifold $Z$.
Using (\FF.2)  one computes that for a submanifold $i:X\ss Z$
$$
i^*\Lambda^{\mu}_q (Z)\ =\ \Lambda^{\mu}_q (X).
\eqno{(\NN.11)}
$$
\medskip

This example  extends directly to the {\bf inhomogeneous} subequation
$$
\l_q(A) \ \geq \ \mu(x)
$$
for a continuous function $\mu(x)$,
by using the local affine jet equivalence $\Phi(A) = A + \mu(x) \cdot I$ to
$ \Lambda^0_q(Z)$. (See Note \LL.11.)

\vskip.3in
%\vfill\eject

%%%%%%%%%%%%%%%%%%%%%%%%%%%%%%%%%%%%%%%%%%%%%%%
%%%%%%%%%%%%%%%%%%%%%%%%%%%%%%%%%%%%%%%%%%%%%%%
%%%%%%%%%%%%%%%%%%%%%%%%%%%%%%%%%%%%%%%%%%%%%%%
%%%%%%%%%%%%%%%%%%%%%%%%%%%%%%%%%%%%%%%%%%%%%%%
%%%%%%%%%%%%%%%%%%%%%%%%%%%%%%%%%%%%%%%%%%%%%%%

\noindent{\headfont Appendix A. Elementary Examples Where Restriction Fails.}\medskip

\bigskip
As noted in Examples \FF.5  and \HH.8 restriction may fail.
Here are two more elementary examples where restriction, and therefore also the Restriction Hypothesis, fail.
In these examples the restricted set  $i^*F$  is closed and hence is a subequation.

\Ex{\KK.1. (First Order)}
   Define $F$ on $\bbr^2$ by $p\pm y^iq^j\geq0$  (where $i$ and $j$ are positive integers).
   Then for the $x$-axis, the restricted subequation $H\equiv i^*F$ is defined by $p\geq0$.

\smallskip
\noindent
{\bf Case $j>i$.}  Restriction to $\{y=0\}$, and hence the restriction hypothesis, fails.
Consider  $u(x,y) = -x + {1\over \a}|y|^\a$ with $\a>0$ small.  Then $p=-1$, and with the
right choice of $\pm$ we have $\pm y^iq^j = |y|^{i+j\a-j}$.  
Thus $p\pm y^iq^j = -1 +|y|^\b\geq0$  with  $\b<0$.
This proves that $u$ is $F$-subharmonic if $|y| >0$ is small.  At points $(x,y) = (x,0)$ there are no
test functions.  Thus $u$ is $F$-subharmonic.  However, the restriction $u\bigr|_X = -x$ is not $H\equiv i^*F$-subharmonic, since $H$ is defined by $p\geq 0$.

\smallskip
\noindent
{\bf Case $j\leq i$.}  The restriction hypothesis, and hence restriction, holds on $\{y=0\}$.
Assume (\EE.5 and 6).  Define $p_\e \equiv p_0 +A_0(x_\e-x_0)$ and $q_\e \equiv {1\over \e}(y_\e-y_0)
= {1\over \e} y_\e$. 
By (\EE.5) we know that $p_\e\pm y_\e^iq_\e^j\geq 0$.
By (\EE.6) we have that $p_\e \to p_0$ and  $|y_\e^iq_\e^j| = {1\over \e^j}|y_\e^{i+j}| \leq  | {y_e^2\over \e}|^j 
\to 0$.  This proves $p_0\geq 0$.\qed

\Ex{\KK.2. (Linear Second-Order and Geometrically Defined)}  
Let $Z \equiv \bbr^{2}$ with coordinates $z=(x,y)$ and
set $X=\{y=0\}$.    Given a section $W(z)$ of  $G(1,\bbr)$ we can write
$W(z) \equiv \span\{\cos \theta(z) e_1 +\sin \theta(z) e_2\}$, defining $\theta(z)$ mod $\pi$.
Then 
$$
P_{W(z)} \ \equiv\ 
\left(
\matrix
{
\cos^2\theta(z)   &   \cos\theta(z) \sin\theta(z)  \cr
 \cos\theta(z) \sin\theta(z)  &  \sin^2\theta(z) 
}
\right).
$$
 The corresponding geometrically defined equation is linear:
 $$
 \bll u \ =\ \tr\left(  D^2 u\bigr|_{W(z)}  \right)
 \ =\ \bra {P_{W(z)}}{  D^2 u_z u}.
 $$
 Set $\sin^2\theta(z) \equiv |y|^\a$.  Then 
 $$
 \bll \ =\ (1-|y|^\a) D_x^2u + 2 |y|^{\a\over 2} (1-|y|^\a)^{1\over 2} D_{x,y}^2 u  +  |y|^\a D_y^2 u.
 $$
 Consider the function 
 $$
 u(x,y) \ \equiv\ -{1\over 2}|x|^2 + {1 \over {2-\b}} |y|^{2-\b}
\eqno{(A.1)}
 $$
 with $0<\a<\b<2$.  At points $y=0$ there are no test functions for $u$.  
 Otherwise $D_x^2 u=-1$, $D_{x,y}^2 u = 0$, and $D_y^2 u = (1-\b)|y|^{-\b}$.
 Hence 
 $$
 \bll u\ = \  -(1-|y|^\a) + {1-\b \over |y|^{\b-\a}}.
 $$
 Since $\a<\b$, $u$ is $\bll$-subharmonic if $|y|$ is small.  However, the restriction satisfies
 $\bll^X \vf = \vf''$,  and $\vf(x) \equiv u(x,0) = -{1\over 2} |x|^2$ is not convex.
 Thus restriction does not hold for $\bll$ even though $\bll$ is linear 
 and $\bll$ is geometrically defined
by the closed subset  $\GG \equiv \{W(x) : z\in \bbr^2\} \ss G(1,\bbr^2)$.
 The restriction hypothesis fails here. Comparing with Theorem 6.4, there
 is no smooth neighborhood retract onto $\GG$;  while comparing with Theorem 5.10, the
 linear restriction hypothesis is satisfied, but the coefficients are not smooth,
 only continuous.

 This counterexample in $\bbr^2$ can be extended to $\bbr^n\times \bbr^m$
 with $u$ still defined by (A.1).  For simplicity, first consider the following linear equation
 even though it is not geometrically defined. The notation is conscripted from (\EE.1)
 $$
\bll \vf \ =\ \tr A + |y|^\a \tr B\ \geq\ 0.
$$
for a constant $\a>0$.  Assume   $\a<\b<2$. 
Note that  $ D^2 ({1\over 2-\b}|y|^{2-\b}) = |y|^{-\b}\{I-\b \hat y\circ \hat y\}$ where $\hat y = y/|y|$.
Hence $\tr \{ D^2 ({1\over 2-\b}|y|^{2-\b})\} = (m-\b)|y|^{-\b}$. For $y\neq0$ we have
 $\bll u = -n + (m-\b)|y|^{\a-\b} \geq 0$. Since $\a-\b<0$,  if $|y|>0$ is sufficiently small,
 then we have $\bll u \geq 0$.
As in $\bbr^2$, $u$ is $\bll$-subharmonic for $|y|$ small, as claimed.

The restricted subequation $H$ on $\{y=0\}$ is just $\Delta_x u\geq0$, which
fails in this case. Hence, restriction and therefore also the restriction hypothesis fail in this case.
We leave it to the reader to find a geometrically defined $\bll$ with
$u$ an $\bll$-subharmonic function.

%\vskip.3in
\vfill\eject

%%%%%%%%%%%%%%%%%%%%%%%%%%%%%%%%%%%%%%%%%%%%%%%
%%%%%%%%%%%%%%%%%%%%%%%%%%%%%%%%%%%%%%%%%%%%%%%
%%%%%%%%%%%%%%%%%%%%%%%%%%%%%%%%%%%%%%%%%%%%%%%
%%%%%%%%%%%%%%%%%%%%%%%%%%%%%%%%%%%%%%%%%%%%%%%
%%%%%%%%%%%%%%%%%%%%%%%%%%%%%%%%%%%%%%%%%%%%%%%

\noindent{\headfont Appendix B. Restriction of Sets of Quadratic Forms Satisfying Positivity.}\medskip

\medskip

In this Appendix  we provide the basic linear algebra material used in our restriction theorems
and their applications.

\bigskip

\centerline{\bf Restriction for Geometrically Determined Subsets of $\Sym(T^*)$}
\medskip

\def\W{W}
\def\H{F}
\def\ff{totally $\H$-free}

Assume that $T$ is an inner product space.  
Let $\Sym(T^*)$ denote the space of quadratic forms on $T$.
Then the {\bf trace} of $A\in\Sym(T^*)$
is well defined, and induces an inner product $\bra A B = {\rm trace}(AB)$ on $\Sym(T^*)$. 
Let $G(p,T)$ denote the grassmannian of $p$-planes in $T$.  By identifying  a subspace a
subspace $V\ss T$ with orthogonal projection $P_V$ onto $V$ we can consider the 
grassmannian $G(p,T)$ to be a subset of  $\Sym(T^*)$.  
Let $i^*A = A\bigr|_V$ denote the restriction of a quadratic form $A\in \Sym(T^*)$ to $V$.
The {\bf $V$-trace} of $A\in \Sym(T^*)$
is defined by 
$$
\tr_V A\ =\ {\rm trace}\left(i^*_V A\right) \ =\ \bra{P_V} A.
$$
\Def{B.1} Given a closed subset $\GG$ of the grassmannian,  the subset $F_\GG \ss  \Sym(T^*)$
defined by 
$$
A\in F_\GG \quad\iff\quad \tr_V A\geq0\ \  \forall\ V\in\GG
\eqno{(B.1)}
$$
is said to be {\bf geometrically determined by $\GG$}.

\medskip

Note that $F_\GG$  is a closed convex cone with vertex at 0.    Moreover, $A\in \Int F_\GG
\iff $ for some $\e>0, \tr_V A\geq \e$ for all $V\in\GG$.   Hence, we have  $F_\GG = \overline{\Int F_\GG}$.
Finally, $F_\GG$ contains no line unless $\GG =\emptyset$, in which case $F_\GG= \Sym(T^*)$.

\Def{B.2}  Given a closed subset  $\GG\ss\G(p,T)$ and a subspace $\W\ss T$ of dimension $\geq p$,
the {\bf $\W$ -tangential part of $\GG$} is defined to be
$$
\GG(\W) \ \equiv\ \{V\in\GG : V\ss \W\}
\eqno{(B.2)}
$$
and we say that $V\in \GG(\W)$ is {\bf tangential to } $\W$.

\Theorem{B.3} {\sl
Suppose that $\H_\GG$ is geometrically determined by the closed subset $\GG\ss G(p,T)$.
Then for each subspace $\W\ss T$ the closure of the restriction of $\H_\GG$ to $\W$ is 
geometrically determined
by the tangential part of $\GG$.  That is
$$
\overline { i^*_\W \H_\GG }  \ =\ \H_{\GG(\W)}.
\eqno{(B.3)}
$$
}

\pf
It suffices to show that 
$$
i^*_W \Int F_\GG\ =\ \Int  F_{\GG(W)}
\eqno{(B.4)}
$$
since $i^*_W F_\GG \ss F_{\GG(W)}$ and $i^*_W \Int F_\GG \ss \Int  F_{\GG(W)}$
are obvious.  (The set $i^*_W \Int F_\GG$ is always open, but $i^*_W F_\GG$ is 
not necessarily closed -- see Example B.6).

Now assume $a\in \Int \H_{\GG(\W)}$.  
Then there exists $\e>0$ such that $\tr_V a \geq \e$ for all $V\in \GG(W)$
Choose $A=\left(\matrix {a& 0\cr0&0\cr}\right)\in\Sym(T^*)$
where the blocking is induced by the splitting $T\equiv \W\oplus N$ with $N=\W^\perp$.  
 Consider the following open neighborhood of $\GG(\W)$ 
in $\GG$
$$
\cn\ \equiv \left\{ V\in\GG : \tr_V A \,>\, \smfrac\e 2\right\}.
\eqno{(B.4)}
$$
Next we use the fact that for all $V\in G(p, T)$
$$
 \bra{P_V}{P_N} \ \geq\ 0\quad{\rm with\ equality\ } \ \ \iff \ \ \ \ V\ss \W.
\eqno{(B.5)}
$$
In particular,
 $$
\inf_{V\in \GG-\cn} \bra {P_V} {P_N}\  \equiv \ \d \ > \ 0.
\eqno{(B.6)}
$$
Set 
 $$
\inf_{V\in \GG-\cn} \bra {P_V} A\  =  \ -M.
\eqno{(B.7)}
$$
Then
 $$
\tr_V(A+tP_N)\ \geq\ -M+t\d \ \ \ \ \ {\rm for} \ \ V\in \GG - \cn
\eqno{(B.8)}
$$
while
 $$
\tr_V(A+tP_N)\ \geq\  \tr_V A \ >\ \smfrac \e 2\ \ \ \ \ {\rm for} \ \ V\in  \cn.
\eqno{(B.9)}
$$
Thus if $t>>0$ so that $-M+t\d >0$, then $A+tP_N \in \Int \H_\GG$, 
and of course $i^*_\W (A - t P_N) = i^*_\W A=a$.\qed

\Def{B.4} The subspace $\W$ is {\bf totally $\GG$-free} if the tangential part of $\GG$ is empty
(i.e., $\GG(\W)=\emptyset$)or equivalently $\H_{\GG(\W)}=\Sym(\W^*)$.
We sat that $F_\GG$ is {\bf unconstrained by} $W$ if $i^*_W F_\GG = \Sym(W^*)$.

\Cor{B.5} {\sl 

\centerline{$\overline{i^*_\W \H_\GG}  = \Sym(\W^*)$}
\smallskip

\centerline{$\iff \ \ i^*_\W \H_\GG = \Sym(\W^*)$\ \ \ \  (i.e.,  $F_\GG$ is unconstrained by  $W$) }
\smallskip
 \centerline{$\iff$\quad
$\GG(\W)=\emptyset$ \ \ \ \ (i.e., $\W$ is totally $\GG$-free).}
}

\pf
Since $\GG(\W)=\emptyset  \iff F_{\GG(W)}= \Sym(W^*)$, it follows from (B.3) that $\overline{i^*_W F_\GG} = \Sym(W^*)
\iff \GG(W)=\emptyset$.  It remains to show that the condition $\overline{i^*_W F_\GG} = \Sym(W^*)$
implies that ${i^*_W F_\GG} = \Sym(W^*)$.  Since $\Int \Sym(W^*) = \Sym(W^*)$, if
$\overline{i^*_W F_\GG} = \Sym(W^*)$, then by (B.4) ${i^*_W F_\GG} = \Sym(W^*)$.\qed

\Ex{B.6. ($i^*_W F_\GG$ is not closed)}
Let $V(s)$ denote the line through $(1,s,s^5)\in\bbr^3$, and $\GG \equiv \{V(s): 0\leq s\leq 1\}$.
Projection onto the line $V(s)$ is given by
$$
P_{V(s)} \ \equiv \ {1\over 1+s^2+s^{10}}
 \left(
\matrix
{
1 & s & s^5  \cr
s & s^2 & s^6  \cr
s^5 & s^6 & s^{10}  \cr
}
\right).
$$
The set $F_\GG$ consists of all $A= ((a_{ij}))\in\Sym(\bbr^3)$ such that
$$
a_{11} +s^2 a_{22} + s^{10} a_{33} +2s a_{12}  +2s^5 a_{13}  +2s^6 a_{23} 
\ \geq\ 0 \fa 0\leq s\leq 1.
\eqno{(B.10)}
$$
Let $W\equiv \bbr^2\times \{0\}$.  Then $\GG(W) = \{V(0)\}$ where $V(0)$ is the line
through $e_1$.  Thus $F_{\GG(W)}$ consists of all
$$
a\ \equiv \left( \matrix 
{
a_{11} & a_{12}  \cr
a_{12} & a_{22}  \cr
}
\right)
\qquad{\rm with}\ \ a_{11}\ \geq\ 0.
$$
In particular, 
$$
a\ \equiv \left( \matrix 
{
0 & 0  \cr
0 & -1  \cr
}
\right)
\ \in\  F_{\GG(W)}.
$$
However, $a\notin i^*_W F_\GG$ because 
$$
A\ \equiv\ 
 \left(
\matrix
{
0 & 0 &  a_{13}  \cr
0 & -1 & a_{23}  \cr
a_{13} &a_{23} & a_{33}  \cr
}
\right).
$$
cannot satisfy (B.10) for small $s>0$.

\vskip .3in

\centerline{\bf Restriction for Subsets of $\Sym(T^*)$ Satisfying Positivity}
\medskip

  Let $\cp\ss\Sym(T^*)$ denote the subset of non-negative quadratic forms.
A subset $F\ss \Sym(T^*)$ is said to satisfy {\bf positivity (P)}  if 
$$
F\,+\,\cp\ \ss\ F.
\eqno{(B.11)}
$$
Of course each $F_\GG$ satisfies (P).

\noindent
\Lemma{ B.7}  {\sl
If   $F\ss \Sym(T^*)$ is a closed set satisfying positivity, then
 \smallskip
 
 (a) \ \ \ \ $F + \Int\cp\ \ss\ \Int F$,
 \smallskip
 
 (b) \ \ \ \ $F = \overline{ \Int F}$,
 \smallskip
 
 (c) \ \ \ \ $\Int F + \cp\ \ss\ \Int F$.
 \smallskip
 \noindent
If, in addition, $F$ is a cone with vertex at the origin, then\smallskip

 (d) \ \ \ \ $F = \Sym(T^*) \quad \iff\quad \exists \, A\in F$ with $A<0$.
 }

\pf 
(a) Note that $A+\Int\cp$ is an open subset of $F$ for each $A\in F$.

(b) \ Pick $P\in \Int \cp$, i.e., $P>0$.  Then by (a) we have that
$A\in F \ \ \Rightarrow\ \ A+\e P\in\Int F$ for each $\e>0$.

(c) Note that $\Int F+P$ is an open subset of $F$ for each $P\in\cp$.

(d) Suppose $F$ contains a negative definite $A<0$.  Then for each 
$B\in\Sym(T^*)$, if $t>>0$ is large enough, $P\equiv B-tA$ is positive. 
Hence, $B= tA+P\in tF+\cp\ss F$.\qed

\Theorem{B.8}  {\sl
Suppose that $\H$ is a closed subset of $\Sym(T^*)$ which is both a cone and satisfies 
(P).  The following conditions on a proper subspace $\W\ss T$ are equivalent.
\smallskip
\noindent
(1) ({\bf $\W$ is $\H$-Morse}) There exists $A\in \H$ with $i^*_\W A<0$.

\smallskip
\noindent
(2) ({\bf $\H$ is unconstrained by $\W$}) $i^*_\W \H =  \Sym(\W^*)$ or equivalently
$\H+\ker i^*_\W =  \Sym(T^*)$.

\smallskip
\noindent
(3)  Given $B\in \ker i^*_\W$, if $B\geq 0$ and rank$\,B = \codim \W$, then$B\in \Int \H$.

\smallskip
\noindent
(3)$'$ ({\bf $\W$ has an $\H$-strict complement}) There exists $B\in \Int \H$ with $i^*_\W B=0$.
}

\Remark{B.9}  If  $\H$ is geometrically defined by $\GG\ss G(p,T)$,
then by Corollary B.5  these conditions are equivalent to the condition that $\W$ contains no $\GG$-planes
($W$ is $\GG$-free).  
 This justifies the following terminology.

\Def{B.10} A subspace $\W$ satisfying the conditions in Theorem B.7 will be called {\bf \ff}.

\pf
 Conditions (1) and (2) are equivalent by (d) above.   Obviously (3) $\Rightarrow$ (3)$'$ since $B\geq0$ with $i^*_\W B=0$ and rank$\,B = \codim \W$ always
  exist.

Next we prove that $(3)'\ \Rightarrow \ (1)$.  If $B\in\Int \H$, then $A\equiv B-\e P \in \H$ with $P>0$ and $\e>0$ small.  If $i^*_\W B=0$, then $i^*_\W A= -\e i^*_\W P <0$ since the restriction of a positive definite quadratic
form is also positive definite.

Finally we show that  $(1)\ \Rightarrow \ (3)$.
Choose $A\in \H$ with $i^*_\W  A < 0$.  Suppose that $B$ satisfies the hypothesis of (3).
Pick $N$ transverse to $\W$ with $T= \W\oplus N$.  Then in block form
$$
A\ \equiv\ \left(  
\matrix
{
-a & c \cr
c^t & b\cr
}
   \right)
   \and 
   B\ \equiv\ \left(  
\matrix
{
0 & \g \cr
\g^t & \b\cr
}
   \right)
$$
where $a= -i^*_\W A>0$ and $0=i^*_\W B$.   Since $B\geq0$, it is a standard fact that $\g=0$.
Since rank$\,B = \dim\,N$, we must have $\b>0$. Set
$$
P\ \equiv\ \smfrac 1 t B - A \ =\ 
\left(  
\matrix
{
a & -c \cr
-c^t & \smfrac 1 t \b -b\cr
}
   \right).
$$
Since $a, \b>0$, one can show that $P>0$ if $t>0$ is sufficiently small.
Hence, $B = tA+tP \in \H +\Int \cp \ss \Int \H$ since $\H$ is a cone
satisfying positivity.\qed
\medskip

Using this algebra one can prove the following topological result which is a vast generalization
of a theorem of Andreotti-Frankel for Stein manifolds.
Given a subequation $F$ on a domain $\O$ we  
define the {\sl free dimension} $\dim_{\rm fr}(F)$
of $F$ to be the largest dimension of a tangent subspace $W\ss T\O$ which is $F$-free.
We say $F$  is   {\sl conical} if each $F_x$ is a cone with vertex at the origin.

\Theorem {B.11} {\sl
Let $F$ be a conical subequation on a domain $\O$ in a manifold $Z$.
 If  $\O$ admits a strictly $F$-subharmonic exhaustion function (i.e.,  if $\O$
is {\bf strictly $F$-convex}), then $\O$ has the homotopy-type of a CW-complex of dimension
$\leq \dim_{\rm fr}(F)$.
}
\pf 
This follows from Morse theory and Theorem B.8 (1) above applied to the
Hessian of the exahustion function at its critical points (cf. [HL$_4$]).\qed

\Remark{B.12}  Let $C^0$ denote the polar of a convex cone $C$.
If $\H\ss \Sym(T^*)$ is a closed convex cone with vertex at the origin
(not necessarily geometrically defined), then for each subspace $\W\ss T$
$$
\H+\ker i^*_\W \ =\ \Sym(T^*) \quad \iff\quad \H^0\cap \Sym(\W^*) \ =\ \{0\},
\eqno{(B.12)}
$$
 since the polar of an intersection is the sum of the polars, and $\ker i^*_\W$ and $\Sym(\W^*)$ are polars 
 of each other.  Thus 
$$
\H^0\cap \Sym(\W^*) \ =\ \{0\} \quad \iff\quad \W\ \ {\rm is\ } \H \ {\rm free}.
\eqno{(B.13)}
$$
This is useful in the convex cone cases which are not geometric.  In the geometric case
$\H_\GG^0 = {\rm ConvexCone}(\GG)$ is the convex cone on $\GG$ with vertex at the origin.
This proves that $\W$ being $\GG$-free can be characterized by either of the following:
$$
\GG \cap  \Sym(\W^*) \ =\ \emptyset \qquad\iff\qquad 
{\rm ConvexCone}(\GG)\cap \Sym(\W)\ =\ \{0\}.
\eqno{(B.14)}
$$

%\vfill\eject
\vskip .3in

%%%%%%%%%%%%%%%%%%%%%%%%%%%%%%%%%%%%%%%%%%%%%
%%%%%%%%%%%%%%%%%%%%%%%%%%%%%%%%%%%%%%%%%%%%%
%%%%%%%%%%%%%%%%%%%%%%%%%%%%%%%%%%%%%%%%%%%%%
%%%%%%%%%%%%%%%%%%%%%%%%%%%%%%%%%%%%%%%%%%%%%

\noindent{\headfont Appendix C.  Extension Results.}

\bigskip

Thus far we have not discussed the extension question:  
\smallskip
%\noindent
{\sl Given a subequation  $F$ on $Z$  and a submanifold $i:X\ss Z$, 
which $i^*F$-subharmonic functions on $X$ are
(locally) the restrictions of $F$-subharmonic functions  on $Z$?}\smallskip

The extreme form of this question arises when $i^*F = J^2(X)$, and so every function
is $i^*F$-subharmonic.  We address this question in two geometrically interesting cases.

Suppose $F\ss J^2(Z)$ is a subequation each fibre of which is a cone
with vertex at the origin ($F$ has the {\sl cone property}). Recall the embedding $\Sym(T^*Z)\ss J^2(Z)$ 
as the 2-jets of functions with critical value zero, and set $F_0 \equiv F\cap \Sym(T^*Z)$.
In Appendix B we have defined what it means for a subspace $W\ss T_zZ$  to be totally
$F_0$-free (see Definition B.5).

\Def{C.1}  A submanifold $X\ss Z$ is said to be {\bf totally $F$-free} if each tangent space
$T_x X$ is totally $F_0$-free.

\Remark{C.2} In the geometric case considered in Section \JJ, a submanifold 
is $F_\GG$-free if it has no tangent $\GG$ planes.
\medskip

In Theorems C.3 and C.6 we assume that $F$ satisfies the mild regularity
condition $\Int (F_x)_0 \ss \Int F$ for each $x\in X$.

\Theorem {C.3}  {\sl
Suppose $F$ is a subequation on $Z$ with  the cone property and that $X\ss Z$ is a closed,  totally 
$F$-free submanifold.  Then every $u\in C^2(X)$ is   the restriction of a strictly
$F$-subharmonic function $\wt u$ on a neighborhood of $X$ in $Z$.
}
\medskip
Now consider a geometric subequation $F_\GG$ on a riemannian $n$-manifold $Z$
determined by $\GG\ss G(p,TZ)$ as in  Section \HH.

\Def{C.4}  A submanifold $X\ss Z$ is strictly $\GG$-convex if at each point $x\in X$
there is a  unit normal vector  $n$ and $\kappa>0$ such that 
$$
\tr_W\left\{  \bra B n \right\}\ \geq\ \kappa \fa W \in \GG(T_{x}X)
\eqno{(C.1)}
$$
where $B$ is the second fundamental form of $X$ (cf. \S \JJ).
(This holds  if $ \GG(T_{x}X)=\emptyset$.)

\Theorem {C.5}  {\sl
Suppose  $X\ss Z$ is a strictly $\GG$-convex submanifold.  
Then every $u\in C^2(X)$ is locally the restriction of a strictly
$\GG$-plurisubharmonic function on $Z$.
}
\medskip

The proof of Theorem C.3 is based on the following result which has other interesting applications.

\Theorem{C.6} {\sl  Suppose that $X$ is a closed submanifold of $Z$, and that $v\in C^2(Z)$ satisfies
\medskip
\centerline{ $X \ =\ \{v=0\}, \qquad v\ \geq\ 0, \qquad$  and\qquad $\rank\,\Hess_x\,v = \codim X, \ \forall \, x\in X$}
\medskip\noindent
Then $X$ is totally $F$-free if and only if the function $v$ is strictly $F$-subharmonic at 
each point of $X$ (and hence in a neighborhood of $X$).
}
\pf
Fix $x\in X$ and set $B\equiv \Hess_x v$.  Then we have 
\medskip
\centerline{ $B\ \geq\ 0, \qquad  B\bigr|_{T_xX}\ =\ 0, \qquad$  and\qquad $\rank\,B\  = \ \codim X$.}
\medskip
If $X$ is totally free, then by Property (3) in Theorem B.2 we have $B\in \Int (F_x)_0$.  Now
since by assumption we have $\Int(F_x)_0 \ss\Int F$, we conclude that $v$ is strictly $F$-subharmonic
at $x$. Conversely, $v$ is strictly $F$-subharmonic  at $x$, then $B \in \Int F_x \cap \Sym(T^*_xZ)$
and $B\bigr|_{T_xX}=0$.  Thus condition (3)$'$ of Theorem B.2 is satisfied, proving that $T_xX$
is $(F_x)_0$-free.\qed

\medskip
\noindent
{\bf Proof of Theorem C.3.}
Pick any $C^2$-extension of  $u$ to $Z$ and also denote it by $u$.
Let $v$ be a function on $Z$ with the properties assumed in Theorem C.6.
We may write $v=\rho^2$ by taking $\rho(z)=\dist(z,X)$ near $X$ for some
riemannian metric on $Z$.  Let $\b:Z\to \bbr$ be a smooth extension
of a given positive function on $X$, and set $\wt u \equiv u+  \b \rho^2$.
  Then we compute that {\sl along the submanifold $X$}:
$$
d \wt u \ =\ du\ 
\and
D^2\wt u \ = \  D^2 u +  \b D^2(\rho^2).
$$ 
That is, along the submanifold $X$:
$$
J(\wt u) = J(u) + \b J( \rho^2).
$$  
At each point $x\in X$ we have $J_x(\rho^2) \in \Int (F_x)_0 \ss \Int F$. 
Therefore by  choosing the positive function $\b$ to be sufficiently large 
at each point $x\in X$, we will have $J(\wt u) \in \Int F$ along $X$, 
and therefore on a neighborhood of $X$ 
in $Z$.\qed

\medskip

\noindent
{\bf Proof of Theorem C.5.}
Fix $x\in X$. It is straightforward to see that by strict $\GG$-convexity,
there is a smooth unit normal vector field $n$ defined in a  compact  neighborhood $V$ of $x$
on $X$ and a $\kappa>0$ so that (C.1) holds at all points of $V$. 

For simplicity we rename $V$ to be $X$.  For clarity we restrict to the case where $Z$ is euclidean space $\rn$ 
Consider the tubular neighborhood
$$
U\ \equiv\ \{  x+\nu \in\rn: x\in X, \nu\in B_\e(0), \nu\perp T_x X\}
$$
for some  small $\e>0$, and define a function $f$ on $U$ by
$$
f(x+\nu) \ =\  \bra {n(x)} \nu + \half c |\nu|^2
$$
 where $c>0$ will be determined later. Set  $\rho(x+\nu)\ =\ \bra {n(x)} \nu$.  Note that $\rho \equiv 0$ on $X$ and therefore
$$
\Hess_X \rho \ \equiv\ 0.
$$
From formula (\JJ.1) we see that 
$$
\Hess_{\rn} \rho \bigr|_{TX} \ =\ \bra B n  \quad {\rm on\ \ }X.
\eqno{(C.2)}
$$
One easily sees that the Hessian of $\half|\nu|^2 = {\rm dist}(\bullet, X)^2$ is 
$$
\half \Hess_{\rn} |\nu|^2 \ =\ P_N\ \equiv\ {\rm orthogonal \ projection\ onto\ the\  normal\ space \ to\ }X
\eqno{(C.3)}
$$
It follows that 
$$
\Hess_{\rn} f\bigr|_{TX} \ =\ \bra Bn.
$$
Hence, by (C.1) we have
$$
\tr_W \left\{   \Hess_{\rn} f \right\} \ \geq\ p\kappa \fa W \in\GG(TX),
$$
and therefore there exists a neighborhood $\cn$ of $\GG(TX)\ss \GG\bigr|_X$ so that
$$
\tr_W \left\{   \Hess_{\rn} f \right\} \ \geq\ \kappa/2  \fa W \in \cn.
$$
Now for a general $W \in \GG\bigr|_X$,
$$
\tr_W \left\{   \Hess_{\rn} f \right\} \ =\ \tr_W \left\{   \Hess_{\rn} \rho \right\} +c \bra {P_W }{P_N}
$$
and by compactness there exists $a>0$ so that
$$
\bra {P_W }{P_N} \ \geq\ a \fa W \in \GG\bigr|_X -\cn.
$$
Let
$$
b\ =\ \inf_{W \in \GG\bigr|_X} \tr_W \left\{   \Hess_{\rn} \rho \right\}.
$$
Then for $c>2|b|/a$ we have
$$
\tr_W \left\{   \Hess_{\rn} f \right\} \ > \ |b| \fa W \in \GG\bigr|_X.
$$
It follows that 
$$
\tr_W \left\{   \Hess_{\rn} f \right\} \ > \ |b| \fa W \in \GG\bigr|_{Nb(X)}
$$
where $Nb(X)$ is a   neighborhood of $X$.  

Now suppose we are given  $u\in C^2(X)$ and $x\in X$. Pick any $C^2$-extension of  $u$ 
to a neighborhood of $X$ and denote it also by $u$.
 On a small compact neighborhood $V$ of $x$ in $X$
apply the construction above to produce the function $f$ on a neighborhood of $V$.
Then for $\l$ sufficiently large, the function $\wt u \equiv  u + \l f$ will be strictly $\GG$-psh on 
a neighborhood of $V$ and satisfy $\wt u\bigr|_V=u$.

For the   case of a general  riemannian manifold $Z$, we use the exponential map
to identify the normal bundle of $X$ with a tubular neighborhood of $X$ in $Z$,
and to the analogous construction.\qed

\vfill\eject

% \magnification=1200
%\NoBlackBoxes
%\nologo 

% \input qtmacros 
% \input QASdefs.tex 

\centerline{\headfont References.}

\vskip .2in

\noindent
\item{[A$_1$]}   S. Alesker,  {\sl  Non-commutative linear algebra and  plurisubharmonic functions  of quaternionic variables}, Bull.  Sci.  Math., {\bf 127} (2003), 1-35. also ArXiv:math.CV/0104209.  

\smallskip

\noindent
\item{[A$_2$]}   S. Alesker,  {\sl  Quaternionic Monge-Amp\`ere equations}, 
J. Geom. Anal., {\bf 13} (2003),  205-238.
 ArXiv:math.CV/0208805.  

\smallskip

\noindent
\item{[AV]}    S. Alesker and M. Verbitsky,  {\sl  Plurisubharmonic functions  on hypercomplex manifolds and HKT-geometry},  J. Geom. Anal. {\bf 16} no. 3 (2006),  375-399.  ArXiv: math.CV/0510140.

\smallskip

\noindent
\item{[Al]}    \ \----------,   {\sl  The Dirichlet problem for the equation Det$\| z_{i,j}\| = \psi(z_1,...,z_n,x_1,...,x_n)$}, I. Vestnik, Leningrad Univ. {\bf 13} No. 1, (1958), 5-24.

\smallskip

\noindent
\item{[BT]}   E. Bedford and B. A. Taylor,  {The Dirichlet problem for a complex Monge-Amp\`ere equation}, 
Inventiones Math.{\bf 37} (1976), no.1, 1-44.

\smallskip

 \item{[B]}  H. J. Bremermann,
    {\sl  On a generalized Dirichlet problem for plurisubharmonic functions and pseudo-convex domains},
          Trans. A. M. S.  {\bf 91}  (1959), 246-276.
\medskip

\noindent
\item{[C]}   M. G. Crandall,  {\sl  Viscosity solutions: a primer},  
pp. 1-43 in ``Viscosity Solutions and Applications''  Ed.'s Dolcetta and Lions, 
SLNM {\bf 1660}, Springer Press, New York, 1997.

 \smallskip

\noindent
\item{[CIL]}   M. G. Crandall, H. Ishii and P. L. Lions {\sl
User's guide to viscosity solutions of second order partial differential equations},  
Bull. Amer. Math. Soc. (N. S.) {\bf 27} (1992), 1-67.

 \smallskip

 \noindent 
\item {[HL$_1$]}   F. R. Harvey and H. B. Lawson, Jr,  {\sl Calibrated geometries}, Acta Mathematica 
{\bf 148} (1982), 47-157.

 \smallskip

\item {[HL$_{2}$]} F. R. Harvey and H. B. Lawson, Jr., 
 {\sl  An introduction to potential theory in calibrated geometry}, Amer. J. Math.  {\bf 131} no. 4 (2009), 893-944.  ArXiv:math.0710.3920.

\smallskip

\item {[HL$_{3}$]} F. R. Harvey and H. B. Lawson, Jr., {\sl  Duality of positive currents and plurisubharmonic functions in calibrated geometry},  Amer. J. Math.    {\bf 131} no. 5 (2009), 1211-1240. ArXiv:math.0710.3921.

\smallskip

\item {[HL$_{4}$]}  F. R. Harvey and H. B. Lawson, Jr., {\sl  Dirichlet duality and the non-linear Dirichlet problem},    Comm. on Pure and Applied Math. {\bf 62} (2009), 396-443.

\smallskip

\item {[HL$_{5}$]} F. R. Harvey and H. B. Lawson, Jr.,  {\sl  Plurisubharmonicity in a general geometric context},  Geometry and Analysis {\bf 1} (2010), 363-401. ArXiv:0804.1316.

\smallskip

\item {[HL$_{6}$]} F. R. Harvey and H. B. Lawson, Jr., {\sl  Dirichlet duality and the non-linear Dirichlet problem
on Riemannian Manifolds},   J. Diff. Geom.  {\bf 88} No. 3 (2011), 395-482.  ArXiv:0907.1981.

\smallskip

\item {[HL$_{7}$]} F. R. Harvey and H. B. Lawson, Jr.,  {\sl  Hyperbolic polynomials and the Dirichlet problem},   ArXiv:0912.5220.
\smallskip

\item {[HL$_{8}$]} F. R. Harvey and H. B. Lawson, Jr., {\sl  Potential theory on almost complex manifolds}, 
  {\sl Ann. Inst.  Fourier} (to appear).  ArXiv:1107.2584.   

\smallskip

\item {[HL$_{9}$]}  \ \----------,   {\sl  Existence, uniqueness and removable singularities for nonlinear
partial differential equations in geometry}, 
Surveys in Geometry (to appear).  
\smallskip

   \noindent
\item{[HM]}    L. R. Hunt and J. J. Murray,    {\sl  $q$-plurisubharmonic functions 
and a generalized Dirichlet problem},    Michigan Math. J.,
 {\bf  25}  (1978),  299-316. 

\smallskip

   \noindent
\item{[K]}    N. V. Krylov,    {\sl  On the general notion of fully nonlinear second-order elliptic equations},    Trans. Amer. Math. Soc. (3)
 {\bf  347}  (1979), 30-34.

\smallskip

   \noindent
\item{[L]}    H. B. Lawson, Jr.,    {Lectures on Minimal Submanifolds}, Vol. I, Second Edition,
 Publish or Perish Press, 
Wilmington, Del., 1980.
\smallskip

\noindent
\item{[NW]}  A. Nijenhuis and W. B. Woolf,
{\sl Some integration problems in almost complex and complex manifolds},    
Ann. of Math. {\bf 77} No. 3 (1963), 424-489.

\smallskip

\item {[P]}  N. Pali, {\sl Fonctions plurisousharmoniques et courants positifs de type (1,1)
sur une vari\'et\'e presque complexe}, 
Manuscripta Math.  {\bf 118} (2005), no. 3, 311-337.

\smallskip

\item {[S]}  Z. Slodkowski, {\sl  The Bremermann-Dirichlet problem for $q$-plurisubharmonic functions},
Ann. Scuola Norm. Sup. Pisa Cl. Sci. (4)  {\bf 11}    (1984),  303-326.

\smallskip

\end

%%%%%%%%%%%%%%%%%%%%%%%%%%%%%%%%%%%%%%%%%%%%%%%%%%
%%%%%%%%%%%%%%%%%%%%%%%%%%%%%%%%%%%%%%%%%%%%%%%%%%
%%%%%%%%%%%%%%%%%%%%%%%%%%%%%%%%%%%%%%%%%%%%%%%%%%
%%%%%%%%%%%%%%%%%%%%%%%%%%%%%%%%%%%%%%%%%%%%%%%%%%
%%%%%%%%%%%%%%%%%%%%%%%%%%%%%%%%%%%%%%%%%%%%%%%%%%
%%%%%%%%%%%%%%%%%%%%%%%%%%%%%%%%%%%%%%%%%%%%%%%%%%
%%%%%%%%%%%%%%%%%%%%%%%%%%%%%%%%%%%%%%%%%%%%%%%%%%
%%%%%%%%%%%%%%%%%%%%%%%%%%%%%%%%%%%%%%%%%%%%%%%%%%

\noindent
\item{[AFS]}  D. Azagra, J. Ferrera and B. Sanz, {\sl Viscosity solutions to second order partial differential
equations on riemannian manifolds}, ArXiv:math.AP/0612742v2,  Feb. 2007.
 \smallskip

\noindent
\item{[BT]}   E. Bedford and B. A. Taylor,  {The Dirichlet problem for a complex Monge-Amp\`ere equation}, 
Inventiones Math.{\bf 37} (1976), no.1, 1-44.

\smallskip

 \item{[B]}  H. J. Bremermann,
    {\sl  On a generalized Dirichlet problem for plurisubharmonic functions and pseudo-convex domains},
          Trans. A. M. S.  {\bf 91}  (1959), 246-276.
\medskip

 \item{[BH]}  R. Bryant and F. R. Harvey,
    {\sl  Submanifolds in hyper-K\"ahler geometry},
          J. Amer. Math. Soc. {\bf 1}  (1989),  1-31.
\medskip

\noindent
 \item{[CNS]}   L. Caffarelli, L. Nirenberg and J. Spruck,  {\sl
The Dirichlet problem for nonlinear second order elliptic equations, III: 
Functions of the eigenvalues of the Hessian},  Acta Math.
  {\bf 155} (1985),   261-301.

 \smallskip

\noindent
\item{[CCH]}  A. Chau, J. Chen and W. He,  {\sl  Lagrangian mean curvature flow for entire Lipschitz graphs},  
ArXiv:0902.3300 Feb, 2009.

 \smallskip

\noindent
\item{[CGG$_1$]}  Y.-G. Chen, Y. Giga and S. Goto  {\sl  Uniqueness and existence of viscosity solutions 
of generalized mean curvature flow equations},  Proc. Japan Acad. Ser. A. Math. Sci  {\bf 65} (1989), 207-210.

 \smallskip

\noindent
\item{[CGG$_2$]}  Y.-G. Chen, Y. Giga and S. Goto  {\sl  Uniqueness and existence of viscosity solutions 
of generalized mean curvature flow equations},  J. Diff. Geom. {\bf 33} (1991), 749-789..

 \smallskip

\noindent
\item{[DK]}   J. Dadok and V. Katz,   {\sl Polar representations}, 
J. Algebra {\bf 92} (1985) no. 2,
504-524.

\smallskip

\noindent
\item{[E]}   L. C. Evans,  {\sl   Regularity for fully nonlinear elliptic equations
and motion by mean curvature},  pp. 98-133 
in ``Viscosity Solutions and Applications''  Ed.'s Dolcetta and Lions, 
SLNM {\bf 1660}, Springer Press, New York, 1997.

 \smallskip

\noindent
\item{[ES$_1$]}   L. C. Evans and J. Spruck,  {\sl   Motion of level sets by mean curvature, I},  
J. Diff. Geom. {\bf 33}  (1991), 635-681.

 \smallskip

\noindent
\item{[ES$_2$]}   L. C. Evans and J. Spruck,  {\sl   Motion of level sets by mean curvature, II},  
Trans. A. M. S.  {\bf 330}  (1992),  321-332.

 \smallskip

\noindent
\item{[ES$_3$]}   L. C. Evans and J. Spruck,  {\sl   Motion of level sets by mean curvature, III},  
J. Geom. Anal.   {\bf 2}  (1992),  121-150.

 \smallskip

\noindent
\item{[ES$_4$]}   L. C. Evans and J. Spruck,  {\sl   Motion of level sets by mean curvature, IV},  
J. Geom. Anal.   {\bf 5}  (1995),   77-114.

 \smallskip

\noindent
\item{[G]}   L. G\aa rding, {\sl  An inequality for hyperbolic polynomials},
 J.  Math.  Mech. {\bf 8}   no. 2 (1959),   957-965.

 \smallskip

\noindent
\item{[Gi]}   Y. Giga, {\sl  Surface Evolution Equations -- A level set approach},
Birkh\"auser,  2006.

 \smallskip

 \noindent 
\item {[H]}   F. R. Harvey,  { Spinors and Calibrations},  Perspectives in Math. vol.9, Academic Press,
Boston, 1990.

 \smallskip

   \noindent
\item{[I]}    H. Ishii,    {\sl  Perron's method for Hamilton-Jacobi equations},    Duke  Math. J.{\bf 55} (1987),  369-384.

\smallskip

\noindent
\item{[PZ]}   S. Peng and D. Zhou, 
{\sl Maximum principle for viscosity solutions on riemannian manifolds},    
ArXiv:0806.4768, June 2008.

\smallskip

\item {[RT]} J. B. Rauch and B. A. Taylor, {\sl  The Dirichlet problem for the 
multidimensional Monge-Amp\`ere equation},
Rocky Mountain J. Math {\bf 7}    (1977), 345-364.

\smallskip

\item {[S]}  Z. Slodkowski, {\sl  The Bremermann-Dirichlet problem for $q$-plurisubharmonic functions},
Ann. Scuola Norm. Sup. Pisa Cl. Sci. (4)  {\bf 11}    (1984),  303-326.

\smallskip

\item {[W]}   J.  B. Walsh,  {\sl Continuity of envelopes of plurisubharmonic functions},
 J. Math. Mech. 
{\bf 18}  (1968-69),   143-148.

\smallskip

\end